\documentclass[dvips,aop]{imsart}

\usepackage{amsmath}
\usepackage[T1]{fontenc}
\usepackage{amsfonts}
\usepackage[english]{babel}
\usepackage[dvips]{graphicx}
\usepackage{graphics}
\usepackage{makeidx}
\usepackage[dvips]{color}
\usepackage{amsmath,amssymb}
\def\wtY{\widetilde Y}
\def\Tr{{\rm Tr}}

\def\wtX{\widetilde X}
\def\wph{{\widetilde p}_h}
\def\wp{\widetilde p}

\def\Dtjjp{\rho^2}

\def\cE{{\mathcal{E}}}

\def\mb0{{\mathbf{0}}}

\def\cA{{\mathcal{A}}}

\def\leftB{[\![}
\def\rightB{]\!]}
\newcommand \A[1]{{\bf (#1)}}

\def\F{{\cal F}}

\def\R{{\mathbb{R}} }
\def\N{{\mathbb{N}} }
\def\E{{\mathbb{E}}  }

\def\P{{\mathbb{P}}  }

\def\det{{\rm{det}}}
\def\tr{{\rm{tr}}}
\def\bint#1^#2{\displaystyle{\int_{#1}^{#2}}}
\def\bsum#1^#2{\displaystyle{\sum_{#1}^{#2}}}

\def\xdt_#1{X_#1(\Delta t)}

\newtheorem{THM}{Theorem}[section]

\newtheorem{PROP}{Proposition}[section]

\newtheorem{LEMME}{Lemma}[section]
\newtheorem{REM}{Remark}[section]

\renewcommand{\thesection}{\arabic{section}}
\newcommand{\mysection}{\setcounter{equation}{0} \section}
\def\finpreuve{
\hfill$\square$
}

\begin{document}

\begin{frontmatter}
\title{ Explicit parametrix and local limit theorems for some degenerate diffusion processes}
\runtitle{Parametrix for some degenerate processes}

\begin{aug}
\author{\fnms{Valentin} \snm{Konakov}\ead[label=e1]{valentin\_konakov@yahoo.com}\thanksref{t1},}
\thankstext{t1}{For the first author, this research was supported by grants 05-01-04004 and 07-01-00077 from the Russian Foundation of Fundamental Researches. This work has partially been written during a visit at the Laboratoire de Probabilit\'es et Mod\`eles Al\'eatoires of the Universities Paris VI and Paris VII in 2007. He is grateful for the hospitality during his stay. Denis Talay is also kindly acknowledged for fruitful discussions.}
\author{\fnms{St\'ephane} \snm{Menozzi} \ead[label=e2]{menozzi@math.jussieu.fr}}
\\\and
\author{\fnms{Stanislav} \snm{Molchanov }
\ead[label=e3]{smolchan@uncc.edu}
}

\runauthor{V. Konakov et Al.}

\affiliation{CEMI RAS, Moscow, Universit\'e Paris VII and
University of North Carolina at Charlotte
}

\address{\printead{e1}\\
Central economics Mathematical institute,\\
Academy of sciences,\\
Nahimovskii av. 47,\\
117418 MOSCOW.\\
RUSSIA.}

\address{\printead{e2}\\
Laboratoire de Probabilit\'es\\
et Mod\`eles Al\'eatoires\\
Universit\'e Paris VII,\\
175 Rue du Chevaleret, 75013 PARIS.\\
FRANCE.
}

\address{\printead{e3}\\
Department of Mathematics\\
University of North Carolina at Charlotte\\
CHARLOTTE.\\
USA.
 }
\end{aug}

\begin{abstract}

For a class of degenerate diffusion processes of rank 2, i.e. when
only Poisson brackets of order one are needed to span the whole space, we obtain a
parametrix representation of Mc Kean-Singer \cite{mcke:sing:67} type for the density. We therefrom derive 
an explicit Gaussian upper bound and a partial lower bound that characterize the
additional singularity induced by the degeneracy.

This particular representation then allows to give a local limit theorem with the usual convergence rate for an associated Markov chain approximation. The key point is that
the "weak" degeneracy allows to exploit the techniques first introduced in
Konakov and Molchanov \cite{kona:molc:84} and then developed in \cite{kona:mamm:00} that rely on Gaussian approximations.

\textit{R\'esum\'e} 

Pour une classe de processus de diffusion de rang deux, i.e. lorsque seuls des crochets de Poisson d'ordre un permettent d'engendrer l'espa\-ce, nous obtenons une repr\'esentation parametrix de type Mc Mean-Singer \cite{mcke:sing:67} de la densit\'e. Nous en d\'erivons une borne sup\'erieure   
Gaussienne explicite et une borne inf\'erieure partielle qui caract\'erisent
la singularit\'e additionnelle induite par la d\'eg\'en\'erescence.

Nous donnons ensuite un th\'eor\`eme limite local pour une approximation par cha\^ine de Markov associ\'ee. Le point crucial est que la faible d\'eg\'en\'erescence permet d'exploiter les techniques initialement introduites par Konakov and Molchanov \cite{kona:molc:84} puis d\'evelopp\'ees dans \cite{kona:mamm:00} et qui reposent sur des approximations Gaussiennes.

\end{abstract}



\end{frontmatter}

\mysection{Introduction}

\subsection{Global overview} \label{sec_11}
Let us consider in $\R^d,\ d\ge 1 $ the Markov diffusion process with generator
$$L=\frac 12 \bsum{i,j\in\leftB 1,d\rightB^2}^{}a_{ij}(x)\partial_{x_ix_j}^2+\bsum{i\in\leftB 1,d \rightB}^{} b_i(x)\partial_{x_i} .$$
If the coefficients of $L$ are smooth enough, say $C^{1}(\R^d) $, bounded, and the diffusion matrix $A(x)=(a_{ij}(x))$ is uniformly elliptic ($ \forall \lambda\in\R^d, \langle A\lambda,\lambda\rangle \in [\delta,\delta^{-1}]$ for an appropriate $\delta>0 $) then the associated process $(X_t)_{t\ge 0}$ has a transition density $p(t,x,y) $ which is the fundamental solution of the parabolic problem
$\partial_t p(.)=L_x p(.), \ p(0,x,y)=\delta_y(x)$. Of course, one also has $\partial_t p(.)=L_y^*p(.),\ p(0,x,y)=\delta_x(y) $.

Moreover, this density satisfies uniformly in $t\in ]0,T]
$ the following Gaussian bounds
$$\frac{M^{-1}}{t^{d/2}}\exp\left(-M\frac{|x-y|^2}{t}\right)\le p(t,x,y)\le \frac{M}{t^{d/2}}\exp\left(-\frac{|x-y|^2}{Mt}\right),$$
where the constant $M$ depends on $T$, $d$, the ellipticity constant and the norms of the coefficients in $C^1(\R^d)$, 
see e.g. Aronson \cite{aron:67} or Stroock \cite{stro:88}.

The above estimations express the following physically obvious fact: if the process starts from
$x_0\in\R^d$, then for small $t>0$, in the neighborhood of $x_0$ it is "almost Gaussian" with the
"frozen" diffusion tensor $A(x_0)$ and the drift $b(x_0)$.

The justification of this fact requires  the solution of the perturbative integral equation for $p(\cdot)$ (so-called \textit{Parametrix equation}), where the leading term of the perturbation theory for $p(\cdot)$ is \textit{exactly} the Gaussian kernel $p_0(\cdot)$ corresponding to the "frozen" coefficients at $x_0$. For details concerning \textit{Parametrix equations} we refer the reader to Mc Kean and Singer \cite{mcke:sing:67}, Friedman \cite{frie:64} or \cite{kona:molc:84}.\\

If the matrix $A(x)$ degenerates, but the coefficients $a,b$ are still smooth, the diffusion process $(X_t)_{t\ge 0}$ with generator $L$ exists (one can use the It\^o calculus for the direct construction of the trajectories), but has generally speaking no density.

Consider now generators of the form $L=\bsum{i=1}^{k}A_i^2+A_0, k<d $, where $(A_i)_{i\in\leftB 0,k\rightB}$ are first order operators (vector fields) on $\R^d$ (or more generally on smooth manifolds) with $C^\infty$ coefficients. Sufficient conditions for the existence of the density can be formulated in terms of the structure of the Lie algebra of the vector fields on $\R^d$, with usual linear operations and the Poisson bracketing $[.,.]$. 
Namely, if the vector fields $A_1,...,A_k,[A_l,A_m]_{(l,m)\in\leftB 0,k\rightB^2}$, $[A_l,[A_m,A_n]]_{(l,m,n)\in\leftB 0,k\rightB^3},...$ span $\R^d$ then the density exists. This result is due to H\"ormander \cite{horm:67}, see also Norris \cite{norr:86} for a Malliavin calculus based probabilistic proof. Operators having the previous property are said to be hypoelliptic. Also, in \cite{horm:67}, H\"ormander stressed that the seed of the idea of hypoellipticity goes back to Kolmogorov's note \cite{kolm:33}.

A. Kolmogorov made the following important observation. Let $d=2$. For the generator $L=\frac 12 \partial_{xx}^2+ax\partial_y,\ a\neq 0$, the solution of the associated SDE writes $(X_t,Y_t)=(x_0+W_t,y_0+a (x_0 t+\int_{0}^{t}W_s ds ))$, where $W$ is a standard one dimensional Brownian motion. Thus $(X_t,Y_t) $ has two dimensional Gaussian distribution with mean $(x_0,y_0+ax_0t)$ and covariance matrix $C=\left(\begin{array}{cc}
t& \frac {at^2}{2}\\
\frac {at^2}{2}& \frac{a^2t^3}{3}
\end{array}
\right)$.
Note that the transition density for small $t$ has higher singularity than the usual heat kernel.
In H\"ormander's form $L=\frac 12 A_1^2+A_0,\ A_1=\partial_x, A_0=a x\partial_y $ so that $[A_1,A_0]=a\partial_y $ and thus, $A_1, [A_1,A_0] $ have together rank 2.

In this paper, using a parametrix approach derived from the work of McKean and Singer \cite{mcke:sing:67}, we are able to derive a Gaussian upper bound, and a "partial" lower bound with the two previous time scales, and an associated local limit theorem
in the following case.

\subsection{Statement of the problem}
We consider $\R^d\times\R^d$-valued diffusion
processes that follow the dynamics
\begin{equation}
\label{DYN_1}
\left\{
\begin{array}{l}
X_t=x+\bint{0}^{t}b(X_s,Y_s) ds+\bint{0}^{t}\sigma(X_s,Y_s) dW_s,\\
Y_t=y+\bint{0}^{t} X_sds,
\end{array}
\right.
\end{equation}
where $(W_t)_{t\ge 0}$ is a standard $d$-dimensional Brownian motion defined on some filtered probability space $(\Omega,\F,(\F_t)_{t\ge 0},\P) $ satisfying the usual assumptions. We assume that $\sigma $ is uniformly elliptic and that $b,\sigma$ are $C^1$, bounded, Lipschitz continuous functions so that there exists a unique strong solution to \eqref{DYN_1}.

Such kind of processes appear in various applicative fields.
For instance in mathematical finance, when dealing with Asian options, $X$ represents the dynamics of the underlying asset and its integral $Y$ is involved in the option Pay-off. Typically, the price of such options write $\E_x[\psi(X_T,T^{-1}Y_T)]$, where for the \textit{put} (resp. \textit{call}) option the function $\psi(x,y)=(x-y)^+ $ (resp. $(y-x)^+ $), see \cite{baru:poli:vesp:01} and \cite{tema:01}. 

The cross dependence of the dynamics of $X$ in $Y$ is also important when handling kinematic models or Hamiltonian systems. 
For a given Hamilton function of the form $H(x,y)=V(y)+\frac{|x|^2}{2}$, where $V$ is a potential and $\frac{|x|^2}{2} $ the kinetic energy of a particle with unit mass, the associated stochastic Hamiltonian system would correspond to $b(X_s,Y_s)= -(\partial_y V(Y_s)+F(X_s,Y_s) X_s)$ in \eqref{DYN_1}, where $F$ is a friction term.
When $F>0$ natural questions arise concerning the asymptotic behavior of $(X_t,Y_t) $, for instance the geometric convergence to equilibrium for the Langevin equation is discussed in Mattingly and Stuart \cite{matt:stua:02}, numerical approximations of the invariant measures in Talay \cite{tala:02}, the case of high degree potential $V$ is investigated in H\'erau and Nier \cite{hera:nier:04}.
Under the previous boundedness assumption on $b$, equation \eqref{DYN_1} describes frictionless Hamiltonian systems with "almost linear" potential.

Importantly, the two time-scales coming from Kolmogorov's example, and that we obtain for the density associated to \eqref{DYN_1}, can be exploited to investigate small time asymptotics of the previous models. For instance, for the Asian option, a normalization is required in the pay-off to make both quantities scale-homogeneous. 

As mentioned above, equation \eqref{DYN_1} provides one of the simplest forms of degenerated processes and the previous assumptions guarantee that H\"ormander's
theorem is satisfied taking only the first Poisson brackets between the vector fields associated to the drift and the diffusive part in \eqref{DYN_1}. 
In a more general hypoelliptic setting, let us mention the work of Cattiaux \cite{catt:90,catt:91} whose assumptions include the case \eqref{DYN_1}, but who obtains less explicit controls, see his Proposition (1.12). Under the "strong" H\"ormander condition that involves the Poisson brackets of the diffusive part of the process, small time asymptotics of the density are discussed in Ben Arous \cite{bena:88} or Ben Arous and L\'eandre \cite{bena:lean:91:2}. Eventually, in whole generality two-sided bounds for the density of degenerate diffusions are investigated in Kusuoka and Stroock \cite{kusu:stro:87}.
All these work strongly rely on Malliavin calculus techniques. We want to stress that the parametrix approach is not very well suited to study general degenerate processes. 
Anyhow, the counterpart is that 
it gives by construction more explicit controls. In the non-degenerate case, for $\alpha$-H\"older continuous coefficients, it directly gives two-sided Gaussian estimates. The lower bound on the diagonal in small time derives from the series representation and the global lower bound is obtained thanks to a chaining argument as in \cite{kusu:stro:87}. Here, we still derive a lower bound in small time from the series, but do not succeed  to do a chaining argument

Also, our controls remain valid if the coefficients in \eqref{DYN_1} are uniformly $\alpha $-H\"older continuous, a case for which H\"ormander's Theorem breaks down, see Section \ref{PREUVE_DIFF} Remark \ref{EXT_NO_HORMANDER} for details.

A natural question then concerns the Markov chain approximation of
\eqref{DYN_1}. For non degenerated processes this aspect has been
widely studied, see e.g. \cite{kona:mamm:00} for local limit
theorems. In \cite{ball:tala:96:2}, using Malliavin calculus
techniques, Bally and Talay
obtain an expansion at order one w.r.t. the time step for the
difference of the densities of the diffusion and a perturbed Euler
scheme, i.e. the stochastic integrals are approximated by Gaussian variables and an artificial viscosity is added to ensure the discrete scheme has a density. This rate corresponds to the usual "weak error" bound. 
Since
we follow the local limit theorem approach we can handle a wider
class of random variables in the approximation but obtain a rate of order
$1/2$ w.r.t the time step. Of course, plugging Gaussian random variables in our approximation yields to rate $h$ as in \cite{ball:tala:96:2}.

Importantly, as opposed to \cite{ball:tala:96:2}, we do not need to introduce an artificial viscosity to ensure the existence of the density for the underlying degenerate Markov chain. We develop analogously to the continuous case a parametrix approach to express the density of the Markov chain in term of the density of an auxiliary frozen random walk. The random walk is degenerated as well, but has a density after a sufficient number of time steps, see Subsection \ref{exi_dens} for details. 
The local limit theorem is then derived from an accurate comparison of the parametrix expansions of the densities of the process and the chain. To motivate this result
we can consider the case of the approximation of a "digital Asian call" i.e. of the quantity $\P[(T^{-1}Y_T-X_T)^+>K]$ for a given $K\in \R^+$. Indeed, the local limit theorem associated to our scheme directly relates the densities of the discrete and continuous objects which is not the case if we only consider a discretization of the non degenerate component and a numerical estimation of the integral, since in  that  case the approximating couple can fail to have a density.

The paper is organized as follows. In Section \ref{ass_not}, we give our assumptions and fix
some notations.
Then, since the form of the Markov chain approximation strongly relies on the proof of our results
for the diffusion we choose to divide this paper into two parts.
Sections \ref{sec_diff} and \ref{PREUVE_DIFF} deal with the results for the diffusion and their proofs.
Section \ref{mark_chain}  is dedicated to the Markov chain approximation of \eqref{DYN_1}, the associated convergence results and the key points of the proofs. The complete proof of the local limit Theorem can be found in the Appendix.  

\subsection{Assumptions and Notations}
\label{ass_not}
We suppose that the coefficients of equation \eqref{DYN_1} satisfy the following assumptions.

\A{UE} $\exists (\lambda_{\min},\lambda_{\max})\in(0,\infty)^2,\
\forall z\in \R^d,\ \forall (x,y)\in\R^{2d},\ \lambda_{\min}|z|^2\le \langle a(x,y) z
,z  \rangle \le \lambda_{\max}|z|^2$, denoting $a(x,y)=\sigma\sigma^*(x,y) $. From now on we suppose that $\sigma $ is the unique symmetric matrix s.t. $\sigma\sigma=a$.  We are interested in the density of the process and its approximation at a given time. Hence, from the uniqueness in law, the previous assumption can be made  without loss of generality.

\A{B} The coefficients $b,\sigma$ in \eqref{DYN_1} are $C^1$, uniformly Lipschitz continuous and bounded.

Throughout the paper we consider the running diffusion \eqref{DYN_1} up to a fixed final time $T>0$. We denote by $C$ a generic positive constant that may change from line to line and only depends on $T$,
and the parameters appearing in \A{UE}, \A{B}.
We reserve the notation $c $ for constants that only depend on parameters from \A{UE}, \A{B}. Other possible dependencies are explicitly indicated.

\mysection{Explicit parametrix and associated controls for the density of the diffusion}
\label{sec_diff}
The previous assumptions guarantee that H\"ormander's Theorem, see e.g. Nualart \cite{nual:98}, holds true, and therefore that $\forall t>0,\ (X_t,Y_t)$ has a density w.r.t. the Lebesgue measure.
Introduce the vector fields in $\R^{2d} $
\begin{eqnarray}
A_0(x,y)=\left(\begin{array}{c}b_1(x,y)\\
\vdots\\
 b_d(x,y)\\
x_1\\
\vdots\\
x_d\\
\end{array} \right),\  \forall j\in\leftB 1,d\rightB,\ A_j(x,y)=\left(\begin{array}{c}
\sigma_{1j}(x,y)\\
\vdots\\
\sigma_{dj}(x,y)\\
0\\
\vdots\\
0\\
\end{array}
\right).\nonumber\\ \label{V_FIELDS}
\end{eqnarray}
One directly derives the following
\begin{PROP}
\label{HORM} For all $(x,y)\in \R^{2d}$,
$${\rm Span}(A_1(x,y),..., A_d(x,y),[A_0(x,y),A_{1}(x,y) ],...,[A_0(x,y),A_{d}(x,y) ])=\R^{2d},$$
where $\forall (i,j)\in \leftB 0,d\rightB^2,\ [A_i,A_j]=A_i \nabla A_j-A_j\nabla A_i $ denotes the Poisson bracket.
\end{PROP}

Fix $T>0$ and $0<t\le T$, $(x,y)\in\R^{2d} $. Since, we
now know that $(X_t,Y_t)$ has a transition density, i.e. $\P[X_t\in
dx',Y_t\in dy'|X_0=x,Y_0=y]=p(t,(x,y),(x',y'))dx'dy'$, our aim is
to develop a parametrix for \eqref{DYN_1} to obtain an explicit
representation of this density.

Recall that we consider the following SDE
\begin{equation}
\label{MODELE_DEV}
\left\{
\begin{array}{l}
dX_s=b(X_s,Y_s)dt+\sigma(X_s,Y_s)dW_s,\ X_0=x,\\
dY_s=X_sds,\ Y_0=y.
\end{array}
\right.
\end{equation}

For the parametrix development we need to introduce a "frozen" diffusion
process, $(\widetilde X_s,\widetilde Y_s)_{s\in [0,t]} $ below. Namely for all $(x',y')\in\R^{2d} $, define
\begin{equation}
\left\{
\begin{array}{l}
d\widetilde{X}_{s}^{t,x',y'}=\sigma (x',y'-x'(t-s))dW_{s}+b(x',y')ds,\widetilde{X}%
_{0}^{t,x',y'}=x, \\
d\widetilde{Y}_{s}^{t,x',y'}=\widetilde{X}_{s}^{t,x',y'} ds, \ \widetilde Y_0^{t,x',y'}=y.
\end{array}
\right.  \label{frozen}
\end{equation}
The key point is that the above process is gaussian. The arguments in the second variable of the diffusion coefficient can seem awkward at first sight, it includes the transport of the frozen point $x'$ with a time reversal. 
This particular choice is actually imposed by the natural metric of the frozen process, see Proposition \ref{CTR_DENS_L}, in order to allow the comparison of the singular parts of the generators. 

The processes $(X_{s},{Y}_{s})$ and $(\widetilde{X}_{s}^{t,x',y'},\widetilde{Y}
_{s}^{t,x',y'}),\ s\in [ 0,t],$ have the following generators: $\forall (x,y)\in\R^{2d},\ \psi\in C^2(\R^{2d} ) $,
\begin{eqnarray}
L\psi(x,y)&=&\biggl(\frac{1}{2}\Tr\left( a(x,y)D_x^2\psi\right) +\langle b(x,y), \nabla_x \psi\rangle 
+\langle x,\nabla _{y}\psi \rangle\biggr)(x,y),\nonumber \\
 \label{DEF_GEN_HAT}
 \widetilde{L}_s^{t,x',y'}\psi(x,y)&=&\biggl(\frac{1}{2}\Tr\left(a\left(x',y'-x'(t-s)\right)D_x^2\psi \right)+\langle b\left(x^{\prime
},y' 
\right)\nabla_x \psi\rangle \nonumber \\
&& 
 +\langle x, \nabla_y\psi\rangle\biggr)(x,y).
 \label{DEF_GEN_TILDE}
\end{eqnarray}

From these operators we define for $0<t\le T, ((x,y),(x',y'))\in(\R^{2d})^2 $ the kernel $H$ by
\begin{equation*}
H(t,(x,y),(x^{\prime },y^{\prime }))=(L-\widetilde{L})\widetilde{
p}(t,(x,y),(x^{\prime },y^{\prime })),
\end{equation*}
where $\widetilde p(t,(x,y),(\cdot,\cdot)):= \widetilde p^{t,x',y'}(t,(x,y),(\cdot,\cdot)),\ \widetilde L:=\widetilde L_{0}^{t,x',y'}$ respectively stand for the density of the process $(\widetilde X_{t}^{t,x',y'},\widetilde Y_{t}^{t,x',y'}) $ and the generator of $(\widetilde X_{s}^{t,x',y'},\widetilde Y_{s}^{t,x',y'})_{s\in[0,t]} $ at time $0$. We omit to explicitly 
emphasize the dependence in $t,x',y'$ for notational convenience.
\begin{REM}
Note carefully that in the above kernel, because of the linear structure of the model the most singular terms, i.e. those involving derivatives w.r.t. $y$, vanish.
\end{REM}

The next proposition gives the expression of the density $ p $
in terms of an infinite sum involving iterated convolutions of the density $\widetilde p $
with the kernel $H$. Namely,
\begin{PROP}[Parametrix expansion for \eqref{MODELE_DEV}]\hspace*{.2cm} \\
\label{PARAM_ANNEXE}
For all $0\le t\le T, ((x,y),(x',y'))\in (\R^{2d})^2 $,
\begin{equation}
\label{EXPR}
p(t,(x,y),(x',y'))=\bsum{r=0}^{+\infty}  \widetilde p\otimes H^{(r)}(t,(x,y),(x',y')),
\end{equation}
where 
\begin{eqnarray*}
f\otimes g(t,(x,y),(x^{\prime },y^{\prime }))&=&\int_{0}^{t}du\int_{\R^{2d}}f(u,(x,y),(z,v))\\
&&\times g(t-u,(z,v),(x^{\prime },y^{\prime }))dzdv,
\end{eqnarray*}
$\tilde p\otimes H^{(0)}=\tilde p$ and $H^{(r)}=H\otimes H^{(r-1)}, \ r>0 $ denotes the $r$-fold convolution of the kernel $H$.
\end{PROP}
The previous Proposition is a direct consequence of the usual parametrix recurrence relations.
For the sake of completeness we provide its proof in Section \ref{PREUVE_DIFF}.

Now, since for $0<t\le T $ $(\tilde X_s,\tilde Y_s)_{s\in[0,t]}$, is a Gaussian process, $\widetilde p $ and its derivatives are well controlled. The previous expression
is the starting point to derive the following
\begin{THM}[Gaussian bounds]\hspace*{.2cm}\\
\label{THEO_DIFF}
There exist constants $c,C>0$ s.t. for all $0\le t\le T, ((x,y),(x',y'))\in (\R^{2d} )^2 $, one has:
\begin{eqnarray} 
p(t,(x,y),(x^{\prime },y^{\prime }))
\leq C 
 \widehat p_{c}(t,(x,y),(x',y'))\nonumber\\
\label{CTR_DENS}
\end{eqnarray}
where $$ \widehat p_{c}(t,(x,y),(x',y')):=\frac{c^{d}3^{d/2}}{(2\pi t^2)^{d}}
\times \exp \left( -c\left[
\frac{\left| x^{\prime }-x\right| ^{2}}{4t}+3\frac{\left| y^{\prime
}-y-\frac{(x+x')t}{2}\right| ^{2}}{t^{3}}\right]\right)$$
enjoys the semigroup property, i.e. $\forall 0\le s< t\le T$,
$$ \int_{\R^{2d}} dw dz \hat p_c(s,(x,y),(w,z))\hat p_c(t-s,(w,z),(x',y'))=\hat p_c(t,(x,y),(x',y')). $$
Also, for a given $C_0>0$, $\exists  t_0:=t_0(C_0,c,C)$ s.t. for $t\le t_0,\ [
\frac{| x'-x| ^{2}}{4t}+3\frac{| y'-y-\frac{(x+x')t}{2}| ^{2}}{t^{3}}]\le C_0,\
p(t,(x,y),(x^{\prime },y^{\prime }))\ge C^{-1} \widehat p_{c^{-1}}(t,(x,y),(x',y'))$.
\end{THM}
\begin{REM}
The lower bound, obtained in small time and compact sets, derives from the parametrix representation of Proposition \ref{PARAM_ANNEXE} and the upper Gaussian control. It remains an open problem to find a well suited chaining argument to derive a global lower bound for this degenerate case.
\end{REM}
\begin{REM}
Note that the above result would remain valid if we replaced the dynamics of $Y_t $ in \eqref{DYN_1} by $Y_t=y+\int_0^t F(X_s) ds$ for a $C^{2+\alpha},\ \alpha>0$, Lipschitz continuous mapping $F:\R^d\rightarrow \R^d$ s.t. $\nabla F\nabla F^* $ is non degenerated, i.e. $\exists c_0, \forall (\xi,x) \in \R^d\times \R^d, |\langle\nabla F\nabla F^*(x) \xi,\xi\rangle|\ge c_0|\xi|^2$. Indeed, in such a case, $(\bar X_s,\bar Y_s)_{s\in[0,T]}:=(F(X_s),Y_s)_{s\in[0,T]} $ follows a dynamics of type \eqref{DYN_1}.
\end{REM}

\mysection{Proof of the main results: diffusion process}
\label{PREUVE_DIFF}

\subsection{Proof of Proposition \ref{PARAM_ANNEXE}: parametrix expansion} 
Following Cattiaux \cite{catt:90} and Lemma \ref{CTR_DENS_L} one derives that $p,\widetilde p $ have continous densities with bounded derivatives.
Hence, from the forward and backward Kolmogorov equations associated to $( X, Y),\ (\widetilde X,\widetilde Y) $ and denoting by $ L^* $ the adjoint of $ L$, we have
\begin{equation*}
p(t,(x,y),(x^{\prime },y^{\prime }))-\widetilde{p}
(t,(x,y),(x^{\prime },y^{\prime }))
\end{equation*}
\begin{equation*}
=\int_{0}^{t}du\frac{\partial }{\partial u}\int_{\R^{2d}}dwdz p
(u,(x,y),(w,z))\widetilde{p}(t-u,(w,z),(x',y'))
\end{equation*}
\begin{equation*}
=\int_{0}^{t}du\int_{\R^{2d} } dwdz\left[ \frac{\partial p
(u,(x,y),(w,z))}{\partial u}\widetilde{p}(t-u,(w,z),(x',y'))
\right.
\end{equation*}
\begin{eqnarray*}
\left. +p(u,(x,y),(w,z))\times \frac{\partial \widetilde{p}(t-u,(w,z),(x^{\prime },y^{\prime
}))}{\partial u}\right]\\
=\int_{0}^{t}du\int_{\R^{2d} }dwdz\left[
{L}^{\ast }{p}(u,(x,y),(w,z))\widetilde{p}(t-u,(w,z),(x^{\prime },y^{\prime }))\right.
\end{eqnarray*}
\begin{equation*}
\left. 
-\widetilde{L
}\widetilde{p}(t-u,(w,z),(x^{\prime },y^{\prime }))p
(u,(x,y),(w,z))\right]
\end{equation*}
\begin{equation*}
=\int_{0}^{t}du\int_{\R^{2d} }dwdz p(u,(x,y),(w,z))({L}-
\widetilde{L})\widetilde{p}(t-u,(w,z),(x^{\prime },y^{\prime }))
\end{equation*}
\begin{equation*}
=p\otimes H(t,(x,y),(x^{\prime },y^{\prime })).
\end{equation*}
A simple iteration completes the proof.\finpreuve
\subsection{Proof of Theorem \ref{THEO_DIFF}}

\subsubsection{Proof of the upper bound}
The proof is divided into two parts. First an elementary control on the density of $(\wtX,\wtY)$ is stated in Lemma \ref{CTR_DENS_L}. Then, this control is used to control the kernel $H$ and the convolutions. 

\textbf{Step 1: Gaussian controls for $(\wtX,\wtY)$.}
\begin{LEMME}
\label{CTR_DENS_L}There exist constants $ c>0,C>0$, s.t. for
all multi-index $\alpha$, $\ |\alpha|\le 3$, $\forall 0
\le u< t \le T $, $\forall (w,z),(x',y')\in \R^{2d}  $
\begin{eqnarray*}
|\partial_w^\alpha \widetilde p(t-u,(w,z),(x',y') )|\leq
C{(t-u)^{-|\alpha|/2}}\frac{c^{d}3^{d/2}}{(2\pi t^2)^{d}}\\
\times \exp\biggl(-c\left[
\frac{|x'-w|^2}{4(t-u)}+3\frac{| y'-z-\frac 12 (x'+w)(t-u) |^2}{(t-u)^3}\right] \biggr)\\
:= C(t-u)^{-\frac{|\alpha|}{2}} \widehat p_c(t-u,(w,z),(x',y')),\\
\widetilde p(t-u,(w,z),(x',y') )\ge 2C^{-1} \widehat p_{c^{-1}}(t-u,(w,z),(x',y')),
\end{eqnarray*}
where $\hat p_c $ is as in Theorem \ref{THEO_DIFF} and enjoys the semi-group property.  
\end{LEMME}
The proof is postponed to Section \ref{LEM_TEC_PREUVE}.

\textbf{Step 2: Control of the kernel.}
Recall that under \A{B}, the coefficients $a,b$ are uniformly Lipschitz continuous. Hence,  it is easy to get from Lemma \ref{CTR_DENS_L} and the previous definition of $H$ that, up to a modification of $c>0$ in $\widehat p_c $,  that  $\exists C_1>0, \forall u\in[0,t) $,
\begin{equation}
|H(t-u,(w,z),(x',y'))|\le \frac{C_1}{\sqrt{t-u}}\widehat p_c(t-u,(w,z),(x',y')).
\label{DER_Y}
\end{equation}
 Lemma \ref{CTR_DENS_L} also yields that $\exists C_2>0,\ \forall u\in (0,t],\ \widetilde p(u,(x,y),(w,z))\le C_2 \widehat p_c(u,(x,y),(w,z))$. 
Setting $C:=C_1\vee C_2 $, we finally obtain
\begin{equation*}
\left| \widetilde{p}\otimes H(t,(x,y),(x^{\prime },y^{\prime
}))\right|
\end{equation*}
\begin{equation*}
\leq \int_{0}^{t}du\int_{\R^{2d} }\widetilde{p}(u,(x,y),(w,z))%
\left| H(t-u,(w,z),(x^{\prime },y^{\prime }))\right| dwdz,
\end{equation*}
\begin{equation*}
\leq \int_{0}^{t}du\int_{\R^{2d} }C^2\widehat p_c(u,(x,y),(w,z))
\frac{1}{\sqrt{t-u}}\widehat p_c(t-u,(w,z),(x',y')) dwdz
\end{equation*}
\begin{equation*}
\leq C^2 t^{1/2} B(1,\frac{1}{2})\widehat p_c(t,(x,y),(x',y')),
\end{equation*}
using the semigroup property of $\widehat p_c $ in the last inequality and where 
$B(m,n)=\int_{0}^{1}du u^{m-1}(1-u)^{n-1} $ denotes the
$\beta$-function. By induction in $r$, %
\begin{eqnarray}
\left| \widetilde{p}\otimes H^{(r)}(t,(x,y),(x^{\prime },y^{\prime
}))\right|
\leq  C^{r+1}t
^{r/2}B(1,\frac{1}{2})B(\frac{3}{2},\frac{1}{2})\times ...\times
B(\frac{r+1}{2},\frac{1}{2})\nonumber\\
\label{THE_EQ_POUR_LES_QUEUES}
\times \widehat p_c(t,(x,y),(x',y')),\ r\in\N^*.
\end{eqnarray}
This implies that the series representing the density $p
(t,(x,y),(x^{\prime },y^{\prime }))$
\begin{equation*}
p(t,(x,y),(x^{\prime },y^{\prime }))=\sum_{r=0}^{\infty }%
\widetilde{p}\otimes H^{(r)}(t,(x,y),(x^{\prime },y^{\prime }))
\end{equation*}
is absolutely convergent and the following estimate holds
\begin{eqnarray*}
\left| p(t,(x,y),(x^{\prime },y^{\prime }))\right|
&\leq& C \widehat p_c(t,(x,y),(x^{\prime },y^{\prime })).
\end{eqnarray*}
 
 \finpreuve
\begin{REM}
\label{EXT_NO_HORMANDER}
Note carefully that  the above series still converges if the coefficients $b,\sigma$ are only uniformly $\alpha$-H\"older continuous. In such case H\"ormander's theorem does not hold, but 
one can show by standard techniques, see e.g. Baldi \cite{bald:78}, that $p(t,(x,y),(.,.)):=\bsum{r\in\N}^{}\widetilde p\otimes H^{(r)}(t,(x,y),(.,.))$ is a probability density and 
derive with a Dynkin like argument, see e.g. Theorem 2.3 in \cite{dynk:63}, that it corresponds to the density of the weak solution of \eqref{DYN_1}. 
\end{REM}

\subsubsection{Proof of the partial lower bound}
\label{LEM_TEC_PREUVE}
From the previous proof and the gaussian nature of $(\widetilde X_t,\widetilde Y_t ) $, see Lemma \ref{CTR_DENS_L}, one gets
\begin{eqnarray*}
p(t,(x,y),(x',y'))&\ge &\widetilde p(t,(x,y),(x',y'))-Ct^{1/2}\widehat p_c(t,(x,y),(x',y'))\\
&\ge& 2C^{-1}\widehat p_{c^{-1}}(t,(x,y),(x',y'))-Ct^{1/2}\widehat p_c(t,(x,y),(x',y'))\\
&\ge& C^{-1}\widehat p_{c^{-1}}(t,(x,y),(x',y'))
\end{eqnarray*}
for $\frac{|x'-x|^2}{4t}+3\frac{| y'-y-\frac 12 (x'+x)t |^2}{t^3}\le C_0 $ and $t$ small enough. 
\subsubsection{Proof of the technical Lemmas} {\phantom{BOUUUUH              }}
\begin{trivlist}
\item[]\textbf{Proof of Lemma \ref{CTR_DENS_L}}.
We prove the statement for $|\alpha|=0$, i.e. without derivation. Indeed, since our computations only involve a finite number of derivations that introduce some polynomials in front of the exponential, the general bound can be derived similarly and the result holds taking the worst constants. Also, with respect to the statement of the lemma, we suppose w.l.o.g. $u=0$ for notational convenience. 
We get from \eqref{frozen} with $x=w,\ y=z $ that for all $0< t\le T $,
\begin{eqnarray}
\wtY_t&=&\left\{z+wt+b\left(x',y'\right)\frac{t^2}{2} \right\}+\bint{0}^{t}\left\{ \int_{0}^{v}\sigma\left(x',y'-x'(t-s) \right)dW_s\right\}dv\nonumber \\
&:=&m_{2,t}+A_{t},\nonumber \\
 A_t&=&\int_0^t (t-s) \sigma\left( x',y'-x'(t-s)\right)dW_s:=\int_0^t (t-s)\widetilde \sigma_{s} dW_s,\nonumber\\\label{expr_loi_YT}
\end{eqnarray}
using It\^o's formula for the last equality. 
Setting $\forall s\in[0,t],\ \widetilde a_s=\widetilde \sigma_s^2$, recall from \A{UE} that $\widetilde \sigma_s$ is symmetric, we finally obtain that  the covariance matrix $\Sigma_t $ of the vector $(\widetilde{X}%
_{t},\widetilde{Y}_{t})$ is equal to
\begin{eqnarray*}
\Sigma_{t}&=&\left(\begin{array}{cc} \int_{0}^t \widetilde a_s ds
&\int_0^t (t-s)\widetilde a_s ds\\
\int_0^t (t-s)\widetilde a_s ds & \int_0^t (t-s)^2\widetilde a_s ds
\end{array}\right).
\end{eqnarray*}
Note from \A{UE} that: $\exists c> 0,\ \forall s\in[0,T],\ \forall \xi\in \R^d,\  \langle \widetilde a_s \xi,\xi \rangle \ge c|\xi|^2$. Hence, by the Cauchy Schwarz inequality
\begin{eqnarray*} 
\forall Z \in \R^{2d}, \langle \Sigma_t Z,Z\rangle \ge c/2 \langle  C_t Z,Z\rangle, \ C_t:=\left(\begin{array}{ll} t I_d& \frac{t^2}{2} I_d\\ \frac{t^2}{2}  I_d& \frac{t^3}{3} I_d  \end{array}\right),
\end{eqnarray*} 
where  $C_t$ is actually the covariance matrix of a $d $-dimensional standard Brownian motion and its integral.
 
The  mean vector  of
$(\widetilde{X}_{t},\widetilde{Y}_{t})$ is equal to $(m_{1,t},m_{2,t})$,
with $  m_{1,t}=w+ b(x',y')t $ and $m_{2,t}$ as in
\eqref{expr_loi_YT}.  
Note that $C _t=T\cA T^*$, where
\begin{eqnarray*} T^{\ast }&=&\left(
\begin{array}{cc }
I_d &\frac{t}{2}I_d\\
0 &I_d
 \end{array}
\right) , \cA=\left( \begin{array}{cc} t I_d &0\\
0 &\frac{t^3}{12}I_d
\end{array}
 \right).
\end{eqnarray*}
Hence,
$C_{t}^{-1}=(T^*)^{-1}\cA^{-1} T^{-1}= \left( \begin{array}{cc} I_d &-\frac t2 I_d\\
0   & I_d
\end{array}\right)  \left(\begin{array}{cc}  t^{-1}I_d& 0\\ 0 & 12 t^{-3}I_d\end{array} \right)\left(
\begin{array}{cc} I_d &0\\
-\frac t2  I_d  &I_d
\end{array}
 \right).$
Now, $\forall Z\in \R^{2d} $,
$\cE:=-\left\langle \Sigma _{t}^{-1}Z,Z\right\rangle \le -c/2\left\langle \cA(T^{-1}Z),T^{-1}Z\right\rangle$.
In particular, for  $
Z=(Z_1,Z_2), Z_1=x'-(w+b(x',y')t),\ Z_2=y'-(z +wt+ b(x',y' )\frac{t^2}{2} )$, we get
 $T^{-1}Z=\left(\begin{array}{c}
x'-w- b(x',y')t\\
 y'-z-\frac 12(x'+w)t 
\end{array}
\right)$.
 We therefore derive
\begin{eqnarray*}
\cE &\le& -\frac{c}{2t}|x'-w-b(x',y')t|^2
-\frac{6c}{t^3}| y'-z-\frac 12(x'+w)t|^2.
\end{eqnarray*}
 From \A{B} (boundedness of $b$), we derive that there exist $c,C>0$ s.t.
 $$\cE \le C- c \left[ \frac{|x'-w|^2}{4t}+3\frac{|y'-z-\frac 12 (w+x')t)|^2}{t^3}  \right]. $$
 Eventually 
 \begin{eqnarray*}
 \widetilde p(t,(w,z),(x',y'))&\le& \frac{Cc^{d}3^{d/2}}{(2\pi t^2)^{d}} \exp\left(- c \left[ \frac{|x'-w|^2}{4t}+3\frac{|y'-z-\frac 12 (w+x')t)|^2}{t^3}  \right]\right)\\
 &:=&C \widehat p_c(t,(w,z),(x',y')).
 \end{eqnarray*}
Note from \cite{kolm:33} that $\widehat p_c $  enjoys the semigroup property.
This gives the statement for $|\alpha|=0 $. The lower bound is derived similarly from the control $\forall Z\in\R^{2d}, \langle \Sigma_t Z,Z\rangle\le \frac{c^{-1}}2\langle C_t Z ,Z \rangle  $ achieved for $c$ small enough.
\finpreuve
  \end{trivlist}
 
\mysection{Markov Chain approximation and associated convergence results}
\label{mark_chain}
\subsection{Global strategy}
Let us recall the strategy to derive a local limit theorem for the Markov chain approximation associated to a diffusion process.
Suppose the underlying diffusion has a density with parametrix representation as in Proposition \ref{PARAM_ANNEXE}. 
If the "natural" Markov chain associated to the diffusion has a density, the main idea is to introduce a Markov chain with frozen coefficients that also has a density so that
the density of the Markov chain can be written in parametrix form as well with a suitable discrete kernel.

The next step consists in comparing these two parametrix representations.
To this end, two key steps are needed: 
\begin{trivlist}
\item[1.] The comparison of the densities of the frozen Markov chain and frozen diffusion process.
\item[2.] The comparison of the kernels. 
\end{trivlist}
The first step relies on Edgeworth like expansions, see e.g. Bhattacharya and Rao \cite{bhat:rao:76}, the second one on careful Taylor like expansions. 

The local limit Theorem is then derived by controlling the iterated convolutions of differences of the kernels. This procedure has been applied successfully in \cite{kona:mamm:00} to derive a local limit theorem for the Markov chain approximation of uniformly elliptic diffusions with bounded coefficients.\\

In our current framework new difficulties arise. First of all it is not obvious to derive that a "natural" Markov chain associated to \eqref{DYN_1} has a density. To guarantee such an existence a common trick 
in the literature consists in adding an artificial viscosity term in the discretization scheme, see e.g. \cite{ball:tala:96:2}. Our strategy is here different. Namely, we manage to obtain a density for the natural frozen Markov chain deriving from \eqref{frozen} after a sufficient number of time steps. We therefore consider a "macro scale" frozen model corresponding to this number of time steps. 
We then obtain a good comparison between the densities of the "aggregated" chain at macro scale and the frozen diffusion process. 
This first step gives the structure of the random variables involved in the approximation in order to have the comparison of the densities. These variables have a density. From these variables, we then derive the Markov chain dynamics by letting the coefficients vary at macro scale.

A second difficulty is that the second component in \eqref{DYN_1} is unbounded. This yields to handle a supplementary term w.r.t. the analysis carried out in \cite{kona:mamm:00} and to a slightly different version of the local limit theorem. In the sequel we first give the dynamics of the Markov chains at macro scale and state the local limit Theorem (Subsection \ref{mod_results}). We give the Lemma for the comparison of the densities (Subsection \ref{comp_dens}) and prove the existence of the density for the aggregated "frozen" Markov chain (Subsection \ref{exi_dens}). The whole proof of the local limit Theorem is carried out in the appendix.

\subsection{Models and results}
\label{mod_results}
Now, fix $T>0,\ \widetilde N\in\N^*$ and let $\tilde h=T/\widetilde N$ be the "micro" time discretization step.
Let $n\in \N^* $ be large enough so that the natural "frozen" chain associated to \eqref{frozen} has a density, see Proposition \ref{EX_DENS}, and define the "macro" scale time step $ h=n\tilde h$ and set $N=\widetilde N/n\in\N^*$ the total number of "macro" time steps over $[0,T] $.

For all $i\in\leftB 0,N\rightB $ set $t_i:=ih$.
For any $(x,y)\in \R^{2d}$, we define on the time grid $\left\{ 0,...,t_N\right\}$ an $\R^{2d} $ valued Markov chain $(Z_{t_i}^h)_{i\in\leftB 0,N\rightB}=((X_{t_i}^h,Y_{t_i}^h)^*)_{i\in\leftB 0,N\rightB}$
whose dynamics is given by
\begin{eqnarray}
Z_{0}^h    &=&(x,y)^{*},\ {\rm and} \ \forall i\in\leftB 0,N-1 \rightB,\ \nonumber \\
X_{t_{i+1}}^h&=& X_{t_i}^h+b(Z_{t_i}^h)h+\sigma(Z_{t_i}^h)\sqrt h\eta_{i+1}^1,\nonumber\\
Y_{t_{i+1}}^h&=& Y_{t_i}^h+(X_{t_i}^h+\frac{\gamma_n}{2}b(Z_{t_i}^h)h+\sigma(Z_{t_i}^h)\sqrt h\eta_{i+1}^2)h,\nonumber\\
\label{CDM}
\end{eqnarray}
where $\gamma_n:=(1+\frac 1n) $.
The variables $(\vartheta_i)_{i \in  (0,N\rightB}:=(\eta^1_{i},\eta_i^2)_{i \in  (0,N\rightB}$ are i.i.d. centered $2d$-dimensional random variables. The density $q_n(\eta_1,\eta_2)$ of $\vartheta_{i} $ satisfies

 \A{A1} $ \E[\vartheta_{i}]=0$, and $Cov(\vartheta_{i})=\left(\begin{array}{ccc}
\mathbf{I}_{d\times d} & \frac 12 \gamma_n \mathbf{I}_{d\times d} \\
\frac 12 \gamma_n \mathbf{I}_{d\times d}& \frac{1}{3}\gamma_n(1+\frac{1}{2n}) \mathbf{I}_{d\times d}
\end{array} \right)
$.

\A{A2} There exist a positive integer $S^{\prime }$ and a function $\psi
:\R^{2d}\rightarrow \R$ with $\sup_{u\in \R^{2d}}\psi (u)<\infty $ \ and $%
\int \left\| u\right\| ^{S}\psi (u)du<\infty $ \ for $S=4dS^{\prime }+4$
such that 
\begin{equation*}
\left| D_{u}^{\nu }q_n(u)\right| \leq \psi (u)
\end{equation*}
for all
$\left| \nu \right|
\in\leftB 0,4\rightB $. 
The main result of the section, i.e. Theorem \ref{THM_LOC_CDM}, is stated in terms of $S'$.



We finally need a "frozen" time inhomogeneous Markov chain.
For
$(x,y),\ (x',y')\in \R^{2d}$,
$j\in  (0,N\rightB^2$ we define
$(\widetilde{Z}^{h}_{t_i})_{i\in\leftB 0, j\rightB}=((\widetilde X_{t_i}^h,\widetilde{Y}_{t_i}^h)^*)_{i\in\leftB 0,j\rightB}$ by
\begin{eqnarray}
\widetilde{Z}_{0}^h&=&(x,y)^{* },\ {\rm and}\ \forall i\in\leftB 0,j-1\rightB\nonumber,\\
\widetilde X_{t_{i+1}}^h&=& \widetilde X_{t_i}^h+b(x',y')h+\sigma(x',y'-x'(t_{j}-t_{i}))\sqrt h\widetilde \eta_{i+1}^1,\nonumber\\
\widetilde Y_{t_{i+1}}^h&=& \widetilde Y_{t_i}^h+ \left\{ \widetilde X_{t_i}^h+\frac{\gamma_n}{2}b(x',y')h+\sigma(x',y'-x'(t_{j}-t_{i}))\sqrt h \widetilde \eta_{i+1}^2\right\}h
.
\nonumber \\
\label{CDM_FRO}
\end{eqnarray}
The i.i.d. variables $(\tilde \eta_i^1,\tilde \eta_i^2)_{i\in (0,j\rightB}
$ have density  $q_n(.)$.

\begin{REM}
Note that the models introduced in \eqref{CDM} and \eqref{CDM_FRO} can seem awkward at first sight. They actually derive from computations that yield the existence of the density for the natural frozen Markov chain associated to \eqref{frozen} after $n$ "micro" time steps $\tilde h$, i.e at the "macro" level with time step $h$. This is developed in Subsection \ref{exi_dens}. 
\end{REM}

From now on,  $p_{h}(t_{j'},(x,y),(\cdot,\cdot))$ and $\widetilde{p}_{h}^{t_{j},x',y'}(t_{j'},(x,y),(\cdot,\cdot)):=\widetilde{p}_{h}(t_{j^{\prime
}},(x,y)$,  $(\cdot,\cdot))$ denote the transition densities between times 0 and $t_{j'}\le t_j$
of the Markov chain \eqref{CDM} and "frozen" Markov chain
\eqref{CDM_FRO} respectively. Introducing a discrete "analogue" to
the generators we derive from the Markov property a relation similar to \eqref{EXPR}
between $ p_h $ and $\widetilde p_h $.

For a sufficiently smooth function $f$, define ${L%
}_{h}$ and $\widetilde{L}_{h\text{ }}$ by
\begin{eqnarray*}
{L}_{h}f(t_{j},(x,y),(x',y'))=\\
h^{-1} \left[ \int
{p}_{h}(h,(x,y),(u,v))f(t_{j}-h,(u,v),(x',y'))dudv\right.\\
 \left.
\phantom{\int_{}^{}} -f(t_{j}-h,(x,y),(x^{\prime
},y^{\prime }))\right] ,
\end{eqnarray*}
\begin{eqnarray*}
\widetilde{L}_{h}f(t_j,(x,y),(x^{\prime },y^{\prime }))=\\h^{-1}%
\left[ \int \widetilde{p}_{h}^{t_{j},x',y'}(h,(x,y),(u,v))f(t_{j}-h,(u,v),(x^{\prime
},y^{\prime }))dudv\right. \\ \left.
\phantom{\int_{}^{}}-f(t_{j}-h,(x,y),(x^{\prime
},y^{\prime }))\right].
\end{eqnarray*}
Note that because of technical reasons, there is a shift in time in the above definitions, i.e. the time is $t_{j}-h$, instead of the "expected"
$t_j$, in the right hand side of the previous equations.

A discrete analogue $H_{h}$ of the kernel $H$ is
defined as
\begin{equation*}
H_{h}(t_j,(x,y),(x',y'))=({L}_{h}-\widetilde{L}_{h})\widetilde{p}_{h}(t_j,(x,y),(x',y')),\ 0<j\le N.
\end{equation*}
From the previous definition
\begin{eqnarray*}
H_{h}(t_j,(x,y),(x^{\prime },y^{\prime }))
=h^{-1}\times \\
\int \left[{p}_{h} -\widetilde{p}%
_{h}^{t_{j},x',y'}\right](h,(x,y),(u,v))
\widetilde{p}_{h}^{t_j,x',y'}(t_{j}-h,(u,v),(x^{\prime
},y^{\prime }))dudv.
\end{eqnarray*}
Analogously to Lemma 3.6 in \cite{kona:mamm:00} we obtain the following result. 

\begin{PROP}[Parametrix 
for
Markov chain].\hspace*{.2cm}\\
Assume \A{UE}, \A{B} are in force. Then, for $0<t_{j}\le T$, 
\begin{equation}
{p}_{h}(t_j,(x,y),(x^{\prime },y^{\prime }))=\sum_{r=0}^{j}\left( \widetilde{p}_{h}\otimes
_{h}H_{h}^{(r)}\right) (t_j,(x,y),(x^{\prime },y^{\prime })),
\label{dev_dens_CDM}
\end{equation}
where the discrete time convolution type operator $\otimes _{h}$ is defined by
\begin{eqnarray*}
(g\otimes _{h}f)(t_{j},(x,y),(x^{\prime },y^{\prime
}))\\ =\sum_{i=0}^{j-1}h\int
g(t_i,(x,y),(u,v))f(t_j-t_i,(u,v),(x^{\prime
},y^{\prime }))dudv,
\end{eqnarray*}
$\widetilde p_h\otimes_h H_h^{(0)}=\tilde p_h$ and $H_h^{(r)}=H_h\otimes_h H_h^{(r-1)}$ denotes the $r$-fold discrete convolution of the kernel $H_h$. W.r.t. the above definition, we use the convention that $\widetilde p_h\otimes_h H_h^{(r)}(0,(x,y),(x',y'))=0, r\ge 1 $.
\end{PROP}
Now \eqref{dev_dens_CDM} and \eqref{EXPR} have the same form.
Comparing these two expressions we obtain the following local limit Theorem.

\begin{THM}[Local limit Theorem for the densities].\hspace{.2cm}\\
\label{THM_LOC_CDM}
Assume \A{UE}, \A{B}, \A{A-1}, \A{A-2} hold true.
Then, $\exists c>0$, 
\begin{eqnarray*}
\sup_{(x,y),(x',y')\in\R^{2d}}\left[(1+|x' |+|x|)\sup_{\delta\in[0,1]}\widehat p_c(T(1+\delta),(x,y),(x',y')) \phantom{\frac {x'+x}{2}}\right.\\
\left. +\chi_{\sqrt{T}}\left(x'-x,y'-y-T\left(\frac {x'+x}{2}\right)\right)\right]^{-1}\\
\times
|p_h(T,(x,y),(x',y')) - p(T,(x,y),(x',y'))| =O(h^{1/2}),
\end{eqnarray*}
where $\widehat p_c $ is as in Theorem \ref{THEO_DIFF}, $p_h $ denotes the density of the Markov chain \eqref{CDM} and $\forall (\rho,u,v)\in\R^+\times \R^{2d} $,
\begin{eqnarray*}
\chi_{\rho}(u,v)=\rho^{-4d}\chi(u/\rho,v/\rho^3),\
 \chi(u,v)=
 \left(1+ (|u|^2+|v|^2)^{
S'-1
} \right)^{-1}
 .
\end{eqnarray*}
\end{THM}
Note from the above result that the bigger is $S'$, the better is the control on the tails.


\subsection{Comparison of the discrete and continuous frozen densities}
 \label{comp_dens}
 The first step for the error analysis is achieved with the following
 \begin{LEMME}\label{THE_LEMME_DENS_GELEE}
There exists $C>0$, s.t. for all $j\in ( 0,N\rightB,\ \rho^2:=t_{j} $,
\begin{eqnarray}
\left| (\wph-\wp) (t_j,(x,y),(x',y')) \right|
\leq Ch^{1/2}\rho ^{-1}\zeta _{\rho }(x^{\prime }-x,y^{\prime }-y-\frac{x+x'}{2}\rho^2),\nonumber\\
\label{DIFF_GELEE}
\end{eqnarray}
where
$\zeta _{\rho }(u,v)=\rho ^{-4d}\zeta (u/\rho ,v/\rho ^{3}),\zeta (u,v)=%
\frac{1}{1+\left[ | u| ^{2}+| v| ^{2}\right] ^{(S-4)/2}%
}$,
$S$ being introduced in \A{A2}.
\end{LEMME}
\textbf{Proof.}
Iterating \eqref{CDM_FRO} from $0$ till
$t_{j} $ we get
\begin{eqnarray}
\label{Dyn_FRO_It}
\widetilde X_{t_{j}}^h
&
=
&
x+b(x',y')\Dtjjp+\rho\{\frac{1}{j^{1/2}}\bsum{k=0}^{j-1}\sigma(x',y'-x'(\rho^2-t_{k}))\widetilde \eta_{k+1}^1 \}\nonumber \\
\widetilde Y_{t_{j}}^h
&
=
&
y+x\Dtjjp+\frac{\rho^4}{2}b(x',y')(1+\frac{1}{nj})\nonumber\\
&&
+\rho^{3}\left\{\frac{1}{j^{1/2}}\bsum{k=0}^{j-1}\sigma(x',y'-x'(\rho^2-t_{k}))\widetilde \eta_{k+1}^2\frac{1}{j}\right. \nonumber \\
&&
+\left.\frac{1}{j^{1/2}}\bsum{k=0}^{j-1}\sigma(x',y'-x'(\rho^2-t_{k}))\widetilde \eta_{k+1}^1\left(1-\frac{k+1}{j}\right) \right\}
\end{eqnarray}
Introduce
\begin{eqnarray*}
m_{j}&=&\left(\begin{array}{cc}
x+ b(x',y')\Dtjjp\\
y+x\Dtjjp+\frac{\rho^4}{2}b(x',y')\gamma_{n,j}
\end{array}\right)
:=
\left(\begin{array}{cc}
m_{j}^1\\
m_{j}^2
\end{array}\right),\ \gamma_{n,j}:=1+\frac{1}{nj},
\end{eqnarray*}
and
\begin{eqnarray*}
\Theta_{j}:=\left(\begin{array}{cc}
\{\frac{1}{j^{1/2}}\bsum{k=0}^{j-1}\sigma(x',y'-x'(\rho^2-t_{k}))\widetilde \eta_{k+1}^1 \}\\
\left\{\frac{1}{j^{1/2}}\bsum{k=0}^{j-1}\sigma(x',y'-x'(\rho^2-t_{k}))\widetilde \eta_{k+1}^2\frac{1}{j}\right. \\
+\left.\frac{1}{j^{1/2}}\bsum{k=0}^{j-1}\sigma(x',y'-x'(\rho^2-t_{k}))\widetilde \eta_{k+1}^1\left(1-\frac{k+1}{j}\right) \right\}
\end{array}\right).
\end{eqnarray*}
The dynamics of \eqref{CDM_FRO} thus writes
$$\left(\begin{array}{c}
\widetilde X_{t_{j}}^h\\
\widetilde Y_{t_{j}}^h
\end{array}\right)=m_{j}+\left(\begin{array}{cc}  \rho\mathbf{I}_{d\times d} & {\mathbf 0}_{d\times d}\\
{\mathbf 0}_{d\times d} & \rho^{3}\mathbf{I}_{d\times d}
\end{array}\right) \Theta_{j}.$$
Setting $\forall s\in[0,\rho^2],\ \phi(s):=\inf\{t_i:=ih: t_i\le s<t_{i+1} \}, \widetilde a_s:=\sigma^2 
(x',y'-x'(\rho^2-s))$ we get  $V_{j}:=Cov(\Theta_{j})=$
\begin{eqnarray*}
\left(\begin{array}{cc} \frac1{t_{j}}\int_0^{t_{j}}ds\widetilde a_{\phi(s)} &  \frac1{t_{j}^2}\int_0^{t_{j}}ds \widetilde a_{\phi(s)}F_1^{j,h} (\phi(s))\\
\frac1{t_{j}^2}\int_0^{t_{j}}ds \widetilde a_{\phi(s)}F_1^{j,h} (\phi(s)) &\frac1{t_{j}^3}\int_0^{t_{j}}ds \widetilde a_{\phi(s)}F_2^{j,h}(\phi(s))
\end{array}\right)
\end{eqnarray*}
where $F_1^{j,h}(\phi(s)):=[\frac{\gamma_n h}2+(t_{j}-(\phi(s)+h ))], F_2^{j,h}(\phi(s)):=[\frac{\gamma_n h^2}3(1+\frac1{2n})+\gamma_n h(t_{j}-(\phi(s)+h ))+(t_{j}-(\phi(s)+h ))^2]$.

 Thus, for $h$ small enough, the covariance matrix $V_{j} $ is uniformly invertible w.r.t. the parameters $n,j,\in\N^* $. 
Denoting by $g_n$  the density of the normalized sum $V_{j}^{-1/2} \Theta_{j}$ we derive
$$\widetilde p_h(t_j,(x,y),(x',y'))=\frac{1}{\rho^{4d}\det(V_{j}^{1/2})} g_n\left(V_{j}^{-1/2}\left(\begin{array}{cc}
\frac{x'-m_{j}^1}{\rho}\\
\frac{y'-m_{j}^2}{\rho^3}
 \end{array}\right)\right).
 $$
Applying the Edgeworth expansion for $g_{n}$ (
the key tool is the normal approximation of
Bhattacharya and Rao, Theorem 19.3 in \cite{bhat:rao:76}) and exploiting arguments similar to those of the proof of Lemma \ref{CTR_DENS_L}, we obtain
\begin{eqnarray}
\left|\widetilde p_h(t_j,(x,y),(x',y'))-\frac{1}{\rho^{4d}\det(V_{j}^{1/2})} g_G\left(V_{j}^{-1/2}\left(\begin{array}{cc}
\frac{x'-m_{j}^1}{\rho}\\
\frac{y'-m_{j}^2}{\rho^3}\end{array}\right)\right)\right|\nonumber \\
\le Ch^{1/2}\rho^{-1}\zeta_\rho(x'-x,y'-y-\frac{x+x'}{2}\rho^2), \label{APRES_ED}
\end{eqnarray}
where $g_G$ stands for the standard $2d$ dimensional Gaussian density. To conclude the proof, recall from the proof of Lemma \ref{CTR_DENS_L} that
\begin{eqnarray}
\widetilde p(t_{j},(x,y),(x',y'))=\frac{1}{\rho^{4d}\det(C_{j}^{1/2})} g_G\left(C_{j}^{-1/2}\left(\begin{array}{cc}
\frac{x'-m_{C,j}^1}{\rho}\\
\frac{y'-m_{C,j}^2}{\rho^3}\end{array}\right)\right)\nonumber\\
\label{COMP_NG_ED}
\end{eqnarray}
where
$m_{C,j}=\left(\begin{array}{cc}
x+ b(x',y')\Dtjjp\\
y+ x\Dtjjp+\frac{\rho^4}{2}b(x',y')
\end{array}\right)
:= \left(\begin{array}{cc}
m_{C,j}^1\\
m_{C,j}^2
\end{array}\right)$,
and $C_{j}=$
\begin{eqnarray*}
\left(\begin{array}{cc}\frac1{t_j}\int_{0}^{t_{j}}ds \widetilde a_s &\frac1{t_j^2}\int_{0}^{t_{j}}ds \widetilde a_s(t_{j}-s) \\ \frac1{t_j^2}\int_{0}^{t_j}ds \widetilde a_s(t_j-s)  & \frac1{t_{j}^3}\int_{0}^{t_{j}}ds \widetilde a_s(t_{j}-s)^2  \end{array}\right).
\end{eqnarray*}
The result eventually follows from \eqref{APRES_ED}, \eqref{COMP_NG_ED} and standard computations involving the mean value theorem.
\finpreuve

\subsection{Existence of the density for the aggregated frozen process}
\label{exi_dens}
Let $h_0>0$ be a given fixed time step. For $i\in \N$ set $t_i:=ih_0 $. Fix $(x',y')\in \R^{2d}, t>0 $.
We consider the frozen model defined by $\widetilde X_0^{h_0}=x,\widetilde Y_0^{h_0}=y $ and for all $i\in\N $,
\begin{eqnarray}
\widetilde{X}_{t_{i+1}}^{h_0} &=&\widetilde{X}_{t_{i}}^{h_0}+b(x',y')h_0+\sigma (x',y'-tx')\sqrt{h_0}\widetilde{\xi }_{i+1},  \nonumber \\
\widetilde{Y}_{t_{i+1}}^{h_0}
&=&\widetilde{Y}_{t_{i}}^{h_0}+\widetilde{X}_{t_{i+1}}^{h_0}h_0
\nonumber \\
&=&\widetilde{Y}_{t_{i}}^{h_0}+h_0
\widetilde{X}_{t_{i}}^{h_0}+h_0^{2}b(x',y')
+h_0^{3/2}\sigma (x',y'-tx')\widetilde{\xi }_{i+1}, \label{FRO_MOD}
\end{eqnarray}
where $(\widetilde \xi_{i})_{i\in\N^*} $ are i.i.d., centered with identity covariance. 
The aim of this section is to show that for $i$ large enough $(\widetilde X_{t_i}^{h_0},\widetilde Y_{t_i}^{h_0}) $ admits a density. We refer the reader to the work of Yurinski \cite{yuri:72} or Molchanov and Varchenko \cite{molc:varc:76} for related topics.

Conditionally to $\left(
\begin{array}{c}
\widetilde{X}_{t_{i}}^{h_0}=x^{\ast } \\
\widetilde{Y}_{t_{i}}^{h_0}=y^{\ast }
\end{array}
\right) $ and iterating the frozen model we get
\begin{eqnarray}
&&\widetilde{X}_{t_{i+n}}^{h_0} =x^{\ast }+(nh_0)b(x',y')+\sigma
(x',y'-x't)\sqrt{nh_0}\widetilde{\xi }_{i,n}^{(1)},\nonumber \\
&&\widetilde{Y}_{t_{i+n}}^{h_0} =y^{\ast }+(nh_0)x^{\ast }+\frac{\gamma_n}{2}(nh_0)^{2}b(x',y')
+(nh_0)^{3/2}\sigma (x',y'-x't)%
\widetilde{\xi }_{i,n}^{(2)},
\label{IT_PETIT}
\end{eqnarray}
where we recall $\gamma_n=(1+\frac{1}{n}) $ and
\begin{eqnarray*}
\widetilde{\xi }_{i,n}^{(1)}&=&\frac{1}{\sqrt{n}
}\left( \widetilde{\xi }_{i+1}+\widetilde{\xi }_{i+2}+...+\widetilde{\xi }
_{i+n}\right),\\
\widetilde{\xi }_{i,n}^{(2)}&=&\frac{1}{\sqrt{n}}\left( \widetilde{\xi }%
_{i+1}+(1-\frac{1}{n})\widetilde{\xi }_{i+2}+...+(1-\frac{n-1}{n})\widetilde{%
\xi }_{i+n}\right) .
\end{eqnarray*}
We have
\begin{eqnarray*}
Var(\widetilde \xi_{i,n}^{(2)})=\frac{(1-\frac{n-1}{n})^{2}+...+1^{2}}{n} &=&\frac{2n^{2}+3n+1}{6n^{2}}=\frac{1}{3}\gamma_n(1+\frac{1}{2n}), \\
Cov(\widetilde \xi_{i,n}^{(1)},\widetilde \xi_{i,n}^{(2)})=\frac{(1-\frac{n-1}{n})+...+1}{n} &=&\frac{n+1}{2n}=\frac{\gamma_n}{2}.
\end{eqnarray*}
Hence, the covariance matrix of the $2d$ dimensional vector $\left(
\widetilde{\xi }_{i,n}^{(1)},\widetilde{\xi }_{i,n}^{(2)}\right) ^{\ast }$
is non-degenerate for $n\ge 2$.

Estimating the characteristic function $\varphi _{n}(\tau _{1},\tau
_{2})$ of the vector $\left( \widetilde{\xi }_{i,n}^{(1)},\widetilde{\xi }%
_{i,n}^{(2)}\right) ^{\ast }\in \R^{2d}$ we derive the following
\begin{PROP}
\label{EX_DENS}
Let $\phi(\tau):=\E\left[\exp\left(i \langle \widetilde \xi_1,\tau \rangle \right)\right],\ \tau\in \R^d$ denote the characteristic function of the $(\widetilde \xi_i)_{i\in\N^*}$. If for all multi index $\beta,\ |\beta|= S+2d+1 $, $|D^\beta\phi(\tau)|\le C(1+|\tau|^{4+2d+1})^{-1} $,
then for $n$ large enough and for all multi index $\alpha $, $|\alpha|\le 4 $, one has
 $$\int_{\R^{d}\times\R^d}^{}|(\tau_1,\tau_2)|^{|\alpha|}|D^{S+2d+1}\varphi_n(\tau_1,\tau_2)|d\tau_1d\tau_2<\infty. $$ In particular, by Fourier inversion the density
\begin{eqnarray}
\label{F_DENS}
\hspace*{1.5cm}f_n(\theta_1,\theta_2)=\frac{1}{(2\pi)^{2d}}\int_{
}^{}\exp(-i\langle (\theta_1,\theta_2)^*,(\tau_1,\tau_2)^*\rangle)\varphi_n(\tau_1,\tau_2)d\tau_1d\tau_2
\end{eqnarray}
 exists and there exists $C$ s.t. for all multi index $\nu,\ |\nu|\le 4 $,
$$|D^\nu f_n(\theta_1,\theta_2)|\le \frac{C}{1+|( \theta_1,\theta_2)|^{S+2d+1}}:=\psi_n(\theta_1,\theta_2).$$
\end{PROP}
\textit{Proof.}
Write
\begin{eqnarray}\varphi _{n}(\tau _{1},\tau _{2}) &=&\E\left[\exp \left\{ i\left\langle \tau _{1},
\widetilde{\xi }_{i,n}^{(1)}\right\rangle +i\left\langle \tau _{2},%
\widetilde{\xi }_{i,n}^{(2)}\right\rangle \right\}\right] 
=\prod_{j=0}^{n-1}\phi \left( \frac{\tau _{1}+(1-\frac{j}{n})\tau _{2}}{%
\sqrt{n}}\right)\nonumber.\\
 \label{PROD_PHI}
\end{eqnarray}
We partition the space $\R^{2d}$ into the following disjoint sets
\begin{eqnarray*}
A_{0} &:= &\left\{ (\tau _{1},\tau _{2})\in \R^{2d}:\left|
\tau _{1}\right| \geq (1-\frac{1}{n})\left| \tau _{2}\right| \right\} , \\
A_{i} &:= &\left\{ (\tau _{1},\tau _{2})\in \R^{2d}:(1-\frac{
i+1}{n})\left| \tau _{2}\right| \leq \left| \tau _{1}\right| <(1-\frac{i%
}{n})\left| \tau _{2}\right| \right\} ,i\in\leftB 1,n-2\rightB, \\
A_{n-1} &:= &\left\{ (\tau _{1},\tau _{2})\in \R^{2d}:\left|
\tau _{1}\right| <\frac{1}{n}\left| \tau _{2}\right| \right\} .
\end{eqnarray*}
If $(\tau _{1},\tau _{2})\in A_{0}$ then for $i\in \leftB 2,n-2\rightB$
\begin{eqnarray*}
\left| \frac{\tau _{1}+(1-\frac{i}{n})\tau _{2}}{\sqrt{n}}\right| &\geq &%
\frac{1}{\sqrt{n}}\left( \left| \tau _{1}\right| -(1-\frac{i}{n})\left|
\tau _{2}\right| \right)  \\
&\geq &\frac{1}{\sqrt{n}}\left( (1-\frac{1}{n})\left| \tau _{2}\right| -(1-%
\frac{i}{n})\left| \tau _{2}\right| \right)   =\frac{i-1}{n\sqrt{n}}\left| \tau _{2}\right|
\end{eqnarray*}
and similarly $\left| \frac{\tau _{1}+(1-\frac{i}{n})\tau _{2}}{\sqrt{n}}\right|  \geq \frac{%
i-1}{n\sqrt{n}}\left| \tau _{1}\right|$. Hence,
\begin{eqnarray}
\label{CTR_A0}
\left| \frac{\tau _{1}+(1-\frac{i}{n})\tau _{2}}{\sqrt{n}}\right|^{2d+1}\ge \frac{(i-1)^{2d+1}}{2n^{3d+3/2}}|(\tau_1,\tau_2)|^{2d+1}.
\end{eqnarray}

If $(\tau _{1},\tau _{2})\in A_{i^{\ast }}$ for some $i^{\ast },$ $
i^{\ast }\in \leftB 1,n-2\rightB$ \ and $l\in\leftB 2,n-1-i^{\ast }\rightB$ then elementary computations yield similarly
\begin{eqnarray}
\label{CTR_AI_1}
\left| \frac{\tau _{1}+(1-\frac{i^{\ast }+l}{n})\tau _{2}}{\sqrt{n}}%
\right| ^{2d+1}
&\geq &\frac{\left( l-1\right) ^{2d+1}}{2n^{3d+3/2}}\left| (\tau _{1},\tau
_{2})\right| ^{2d+1},
\end{eqnarray}
and for $l\in\leftB 1,i^{\ast }-1\rightB$
\begin{eqnarray}
\label{CTR_AI_2}
\left| \frac{\tau _{1}+(1-\frac{i^{\ast }-l}{n})\tau _{2}}{\sqrt{n}}%
\right| ^{2d+1} &\geq& \frac{l^{2d+1}}{2n^{3d+3/2}}\left| (\tau _{1},\tau _{2})\right| ^{2d+1}.
\end{eqnarray}

If \ $(\tau _{1},\tau _{2})\in A_{n-1}$ then for $i\in\leftB 1,n-1\rightB$
\begin{eqnarray}
\label{CTR_AN}
\hspace*{.95cm}\left| \frac{\tau _{1}+(1-\frac{i}{n})\tau _{2}}{\sqrt{n}}\right| ^{2d+1}
&\geq &\frac{1}{2n^{d+1/2}}\left( 1-\frac{i+1}{n}\right) ^{2d+1}\left| (\tau _{1},\tau
_{2})\right| ^{2d+1}.
\end{eqnarray}

Use now the growth assumption on $\phi$ and the inequality $1+\sum_{j=1}^{N}p_{j}\leq \prod_{j=1}^{N}(1+p_{j})$
where $p_{j}\geq 0,$ to derive from \eqref{PROD_PHI}
\begin{eqnarray*}
\left| \varphi _{n}(\tau _{1},\tau _{2})\right|  &=&\left|
\prod_{j=0}^{n-1}\phi \left( \frac{\tau _{1}+(1-\frac{j}{n})\tau _{2}}{\sqrt{%
n}}\right) \right| \leq \frac{C^n}{\prod_{j=0}^{n-1}\left( 1+\left| \frac{%
\tau _{1}+(1-\frac{j}{n})\tau _{2}}{\sqrt{n}}\right| ^{2d+1}\right) } \\
&\leq &\frac{C^n}{1+\sum_{j=0}^{n-1}\left| \frac{\tau _{1}+(1-\frac{j}{n}%
)\tau _{2}}{\sqrt{n}}\right| ^{2d+1}}.
\end{eqnarray*}
Now equations \eqref{CTR_A0}, \eqref{CTR_AI_1}, \eqref{CTR_AI_2}, \eqref{CTR_AN} yield that there exists $n$ large enough s.t.
\begin{eqnarray*}
\left| \varphi _{n}(\tau _{1},\tau _{2})\right| \leq \frac{C(n)}{1+\left| (\tau _{1},\tau
_{2})\right| ^{2d+1}},
\end{eqnarray*}
where $C(n)\underset{n}{\rightarrow} +\infty$.
Anyhow, for such a fixed $n$, one has $\varphi_n \in L^1(\R^{2d})$ which
implies the existence of the density $f_n$ of the vectors  $\left( \widetilde{\xi }_{i,n}^{(1)},\widetilde{%
\xi }_{i,n}^{(2)}\right) ^{\ast }\in \R^{2d}$. The properties concerning the growth and derivatives of $f_n$ are derived from \eqref{F_DENS} and the growth and smoothness properties of $\phi $. \finpreuve

Hence we can set $(\eta_i^1,\eta_i^2):=(\widetilde\xi_{i,n}^{(1)},\widetilde\xi_{i,n}^{(2)}) $ where $(\widetilde\xi_{i,n}^{(1)},\widetilde\xi_{i,n}^{(2)})$ are as in the above proposition. Introducing  a "macro" scale time step $h=nh_0$, the discrete model \eqref{CDM_FRO} corresponds to the "aggregated" dynamics of \eqref{IT_PETIT}.  
Set for all $(\theta_1,\theta_2)\in \R^{2d},\ \psi(\theta_1,\theta_2):=\psi_n(\theta_1,\theta_2)$. With the notations of Section \ref{mod_results} one derives that $q_n(\theta_1,\theta_2)$ $=f_n(\theta_1,\theta_2)$ satisfies \A{A2} with the above $\psi $.

\appendix
\section{Proof of the local limit theorem \ref{THM_LOC_CDM}}
From now on, we use the following
notations for multi-indices and powers. For
$\nu =(\nu _{1},...,\nu _{2d})\in \N^{2d},$ $%
(x,y)=(x_{1},...,x_{d},y_{1},...,y_{d})^{* }$ set
\begin{eqnarray*}
&&| \nu | =\nu_{1}+... +\nu _{2d},\ \nu!=\nu_1!...\nu_{2d}!, \\
&& \text{ }(x,y)^{\nu
}=x_{1}^{\nu _{1}} ...\ x_{d}^{\nu_{d}}y_1^{\nu_{d+1}} ...\ y_d^{\nu_{2d}},
D^{\nu }=D_{x_{1}}^{\nu
_{1}} ... D_{x_{d}}^{\nu _{d}}D_{y_1}^{\nu_{d+1}}...D_{y_d}^{\nu_{2d}}. 
\end{eqnarray*}

\subsection{Preliminary controls on the discrete kernel}
We first give some controls for the kernel $H_{h}(t_{j},(x,y),(x^{\prime
},y^{\prime }))$. Namely, the following Lemma states that the difference between $H_h$, $K_h:=(L-\widetilde L )\widetilde p_h$ and an additional remainder term $M_h$ is small, i.e. has the order announced in Theorem \ref{THM_LOC_CDM}.
\begin{LEMME}[Control of the discrete kernel]
\label{CTR_KERNL_MARK}
For all $j\in\leftB 1,N\rightB $, set $\rho^2=t_j $. One has
\begin{eqnarray}
\left| (H_{h}-K_{h}-M_{h})(t_j,(x,y),(x^{\prime },y^{\prime }))\right|\nonumber\\
\leq Ch^{1/2}\rho ^{-1}\zeta _{\rho }(x^{\prime }-x,y^{\prime }-y-\frac{x+x'}{2}\rho^2).
\label{anc_10}
\end{eqnarray}
where $\zeta_\rho$ is as in Lemma \ref{THE_LEMME_DENS_GELEE} and for $j\in(1,N\rightB$,
\begin{equation*}
K_{h}(t_j,(x,y),(x^{\prime },y^{\prime }))=({L}-%
\widetilde{L})\widetilde{p}_{h}(t_j,(x,y+xh),(x^{\prime
},y^{\prime })),
\end{equation*}
i.e. $K_h$ is the difference of the generators associated to the initial and frozen
diffusion processes between $0 $ and $t_{j} $ applied to the density of the Markov chain with a slight change for the initial point in the $y$ component, 
\begin{equation}
M_{h}(t_j,(x,y),(x^{\prime },y^{\prime
}))=\bsum{k=1}^{3}M_{h}^k(t_j,(x,y),(x^{\prime
},y^{\prime })), \label{anc_11}
\end{equation}
where the $(M_{h}^k)_{k\in\leftB 1,3\rightB} $ are defined in 
Appendix \ref{APPENDICE}.

For $j=1$ we set $K_{h}(t_j,(x,y),(x^{\prime },y^{\prime }))=0$,
\begin{equation*}
M_{h}(t_j,(x,y),(x^{\prime },y^{\prime }))=H_h(t_j,(x,y),(x^{\prime },y^{\prime })).
\end{equation*}
\end{LEMME}
The proof is postponed to Appendix \ref{APPENDICE}. From this proof one also derives that the terms appearing in Lemma \ref{CTR_KERNL_MARK} are controlled with the following:
\begin{LEMME}
\label{OLD_LEMME2}
There exists a constant $C$ s.t. for all
$0<j\le N$, for all $(x,y)$ and $(x^{\prime
},y^{\prime })$ in $\R^{2d} $
\begin{eqnarray*}
 (|K_{h}|+|M_h|+\bsum{i=1}^{3}|M_h^i|+|H_h|)(t_j,(x,y),(x^{\prime },y^{\prime })\\
\leq C
\rho ^{-1}
\zeta _{\rho }\left( x^{\prime }-x,y^{\prime
}-y-\frac{x+x'}{2}\rho^2\right) ,
\end{eqnarray*}
with $\zeta _{\rho }$ as in Lemma \ref{THE_LEMME_DENS_GELEE}. Here again $\rho =\sqrt{t_j}.$
\end{LEMME}
The key fact is that the previous bound provides an integrable singularity in $\rho$.

\subsection{Comparison of the parametrix expansions for the compensated diffusion and Markov chain}
We first state an auxiliary result concerning the behavior of the iterated discrete kernel applied to the density of the frozen Markov chain.
\begin{LEMME}
\label{OLD_LEMME3}
There exists a constant $C$ (that does not depend on $(x,y)$ and $(x^{\prime
},y^{\prime })$) such that, for all $0<j\le N,\ r\in\leftB 0,j\rightB $,
\begin{eqnarray*}
\left| \left( \widetilde{p}_{h}\otimes _{h}H_{h}^{(r)}\right)
(t_j,(x,y),(x^{\prime },y^{\prime }))\right| &\leq & \frac{C^{r+1}\text{ }\rho ^{r}%
}{\Gamma \left( 1+\frac{r}{2}\right) }\\
&& \times \chi _{\rho }\left( x^{\prime
}-x,y^{\prime }-y-\frac{x+x'}{2}\rho^2\right),
\end{eqnarray*}
where 
$\chi_\rho $ and $S'$ are as in Theorem \ref{THM_LOC_CDM}.
\end{LEMME}
To prove the lemma it is sufficient to repeat the proof of
Lemma 3.11 in \cite{kona:mamm:00} with obvious modifications
concerning the additional arguments in $y,y^{\prime }$.

\begin{LEMME}
\label{OLD_LEMME4}
For $0<j\leq N$
the following formula holds:
\begin{equation*}
p_{h}(t_j,(x,y),(x',y'))=\sum_{r=0}^{j}\left( \widetilde{p}\otimes
_{h}(M_{h}+K_{h})^{(r)}\right) (t_j,(x,y),(x^{\prime
},y^{\prime }))+R,
\end{equation*}
where
$\left| R\right| \leq Ch^{1/2}\rho ^{-1}\chi _{\rho }(x^{\prime
}-x,y^{\prime }-y-\frac{x+x'}{2}\rho^2)$
for some constant $C$. The function $\chi _{\rho }$
is as in Theorem \ref{THM_LOC_CDM}.
\end{LEMME}

The proof follows from Lemmas \ref{THE_LEMME_DENS_GELEE} and \ref{OLD_LEMME2} and is analogous to the proof
of Lemma 3.13. in \cite{kona:mamm:00}.\finpreuve

Let us now compare the parametrix expansions of the compensated diffusion and Markov chain at time $T$.  From Proposition \ref{PARAM_ANNEXE}, \eqref{THE_EQ_POUR_LES_QUEUES} and Stirling's asymptotic formula for the $\Gamma$ function
we have
\begin{equation}
p(T,(x,y),(x',y'))=\sum_{r=0}^{N}\left( \widetilde{p}\otimes H^{(r)}\right)
(T,(x,y),(x^{\prime },y^{\prime }))+R_{1},  \label{anc_21}
\end{equation}
where $\left| R_{1}\right| \leq Ch^{1/2}\widehat p_{c}(T,x'-x,y'-y-\frac{x+x'}{2} T)$,
 with $\widehat p_c  $ as in Theorem \ref{THEO_DIFF}. By Lemma \ref{OLD_LEMME4}
\begin{equation}
p_{h}(T,(x,y),(x^{\prime },y^{\prime }))=\sum_{r=0}^{N}\left( \widetilde{p}%
\otimes _{h}(M_{h}+K_{h})^{(r)}\right) (T,(x,y),(x^{\prime
},y^{\prime }))+R_{2}  \label{anc_22}
\end{equation}
where
\begin{equation*}
\left| R_{2}\right| \leq Ch^{1/2}T^{-1/2}\chi _{\sqrt{T}}(x^{\prime
}-x,y^{\prime }-y-\frac{x+x'}{2}T).
\end{equation*}
Because of \eqref{anc_21} and \eqref{anc_22}, to prove the theorem it remains
to show that
\begin{eqnarray}
|\Delta_N|&:=&\left| \left\{\sum_{r=0}^{N}\left( \widetilde{p}\otimes H^{(r)}\right)
-
\sum_{r=0}^{N}\left( \widetilde{p}%
\otimes _{h}(M_{h}+K_{h})^{(r)}\right)\right\}(T,(x,y),(x^{\prime
},y^{\prime }))  \right|\nonumber\\
&&\leq C(1+|x'|)h^{1/2}\chi _{\sqrt{T}}(x^{\prime }-x,y^{\prime }-y-\frac{x+x'}{2}T).
\label{anc_23}
\end{eqnarray}
Note that $|\Delta_N|\le S_1+S_2+S_3+S_4 $, where
\begin{equation*}
S_{1}=\left| \left(\sum_{r=0}^{N}\left( \widetilde{p}\otimes
H^{(r)}\right)
-\sum_{r=0}^{N}\left( \widetilde{p}%
\otimes _{h}H^{(r)}\right)\right) (T,(x,y),(x^{\prime },y^{\prime
}))\right| ,
\end{equation*}
\begin{equation*}
S_2=\left| \left(\sum_{r=0}^{N}\left( \widetilde{p}\otimes_h
H^{(r)}\right)
-\sum_{r=0}^{N}\left( \widetilde{p}%
\otimes _{h}\widetilde H^{(r)}\right)\right) (T,(x,y),(x^{\prime },y^{\prime
}))\right|,
\end{equation*}
\begin{equation*}
S_{3}=\left| \left( \sum_{r=0}^{N}\left( \widetilde{p}\otimes
_{h}\widetilde H^{(r)}\right)-\sum_{r=0}^{N}\left( \widetilde{p}%
\otimes _{h}(M_{h}+\widetilde H)^{(r)}\right)\right) (T,(x,y),(x^{\prime },y^{\prime
}))\right| ,
\end{equation*}
\begin{eqnarray*}
S_{4}&=&\left|\left(\sum_{r=0}^{N}\left( \widetilde{p}\otimes
_{h}(M_{h}+\widetilde H\right)
^{(r)}
-\sum_{r=0}^{N}\left( \widetilde{p%
}\otimes _{h}(M_{h}+K_{h})^{(r)}\right)\right) (T,(x,y),(x^{\prime
},y^{\prime }))\right|,
\end{eqnarray*}
where $\tilde H(t,(x,y),(x',y'))=H(t,(x,y+xh),(x',y')) $ is a shifted operator introduced
for the comparison with $K_h$, see the proof of Lemma \ref{CTR_KERNL_MARK} in the Appendix \ref{APPENDICE} for details.

We shall show
\begin{eqnarray*}
S_{i}&\leq&  Ch^{1/2}\chi _{\sqrt{T}}(x^{\prime }-x,y^{\prime }-y-\frac{x+x'}{2}T),\
i\in \{1,3,4 \}, \\
S_2&\le &C(1+|x'|)h^{1/2}\widehat p_c (T,(x,y),(x^{\prime },y^{\prime })).
\end{eqnarray*}
This is done in Appendix \ref{CTR_SI}.

\section{Proof of Lemmas \ref{CTR_KERNL_MARK} and \ref{OLD_LEMME2}}
\label{APPENDICE}

\subsection{Proof of Lemma \ref{CTR_KERNL_MARK}}
For $j=1$ we have $\rho =\sqrt{h}.$ By definition of $H_{h}$
\begin{equation*}
H_{h}(h,(x,y),(x',y'))=\rho ^{-2}\left[ ({p}
_{h}-\widetilde{p}_{h}^{h,x',y'})(h,(x,y),(x^{\prime },y^{\prime
}))\right] .
\end{equation*}
Thus, recalling $q_n$ stands for the density of the random variables appearing in schemes \eqref{CDM}, \eqref{CDM_FRO}
\begin{eqnarray*}
\left| H_{h}(h,(x,y),(x^{\prime },y^{\prime }))\right|
=h^{-(1+2d)}\left| \frac{1}{\sqrt{\det
a(x,y)}}q_{n}\left( u+\delta _{1},v+\delta_2\right)
\right. 
\end{eqnarray*}
\begin{equation*}
\left.-\frac{1}{\sqrt{\det
a(x',y'-x'h) } }  q_{n}\left( u,v\right)\right|, 
\end{equation*}
where 
\begin{equation*}
u=\frac{\sigma ^{-1}(x',y'-x'h)(x'-x-b(x',y')h)}{\sqrt{h}}%
,u+\delta _{1}=\frac{\sigma ^{-1}(x,y)(x'-x-b(x,y)h)}{\sqrt{h}},
\end{equation*}
\begin{eqnarray*}
v&=&\sigma^{-1}(x',y'-x'h)\frac{y^{\prime }-y-(x+\frac{1}{2}\gamma _{n}hb(x',y'))h}{h^{3/2  }},\\
v+\delta _{2}&=&
\sigma^{-1}(x,y)\frac{y^{\prime }-y-(x+\frac{1}{2}\gamma _{n}hb(x,y) )h}{h^{3/2 }}.
\end{eqnarray*}
Note that 
\begin{equation}
\left| \frac{1}{\sqrt{\det a(x',y'-x'h)}}-\frac{1}{\sqrt{\det a(x,y)}}\right| \leq
C\left[h^{1/2}\left( \left| u\right| +h^{1/2}\right)+h^{3/2}\left( \left| v\right| +h^{1/2}\right)\right],  \label{OLD_1A}
\end{equation}
\begin{equation}
\left| q_{n}\left( u+\delta _{1},v+\delta_2\right)   -q_{n}\left( u,v \right)  \right|
\leq  C\psi \left( \left[ u,u+\delta _{1}\right] _{\gamma },\left[
v,v+\delta _{2}\right] _{\gamma }\right) \left( \left| \delta _{1}\right|
+\left| \delta _{2}\right| \right),   \label{OLD_1B}
\end{equation}
where $\left[ X,Y\right] _{\gamma }:= (1-\gamma )X+\gamma Y,$ $%
\gamma \in \lbrack 0,1],$ for $X,Y$ matrices or vectors. One has
\begin{eqnarray*}
\left[ u,u+\delta _{1}\right] _{\gamma }&=&\frac{[\sigma ^{-1}(x',y'-x'h),\sigma ^{-1}(x,y)]_{\gamma }(x'-x-b(x',y')h)}{\sqrt{h}}\\
&&+\gamma \sigma ^{-1}(x,y)(b(x',y')-b(x,y))\sqrt{h}\nonumber\\
&:=& \lbrack \sigma ^{-1}(x',y'-x'h),\sigma ^{-1}(x,y)]_{\gamma }\sigma
(x',y'-x'h)u+\gamma R_1,  
\end{eqnarray*}
\begin{eqnarray*}
\left[ v,v+\delta _{2}\right] _{\gamma }&=&\frac{[\sigma ^{-1}(x',y'-x'h),\sigma ^{-1}(x,y)]_{\gamma }(y'-y-(x+\frac{\gamma_nb(x',y')h}2)h)}{h^{3/2}}\\
&&+\gamma \sigma ^{-1}(x,y)(b(x',y')-b(x,y))\frac{\gamma_n}2\sqrt{h}\nonumber\\
&:=& \lbrack \sigma ^{-1}(x',y'-x'h),\sigma ^{-1}(x,y)]_{\gamma }\sigma
(x',y'-x'h)v+\gamma R_2.  
\end{eqnarray*}
Assumptions \A{B}, \A{UE} also yield
\begin{equation*}
| R_1|+|R_2 |\le Ch^{1/2}, 
\end{equation*}
\vspace*{-0.8cm}
\begin{eqnarray*}
&&\exists C_1,C_2>0,\ \forall \xi\in\R^d,\ 0\leq \gamma \leq 1,\\
&& C_{1}| \xi| \leq | \lbrack \sigma ^{-1}(x',y'-x'h),\sigma
^{-1}(x,y)]_{\gamma }\sigma (x',y'-x'h)\xi| \leq C_{2}| \xi|.  
\end{eqnarray*}
We also have 
 \begin{equation*}
| \delta _1| +| \delta _2| \leq Ch^{\frac{1}{2}}\left( 1+| u|^2+|v  |^2\right) .  
\end{equation*}
Recall that from \A{A2} and our choice for $\psi$, 
$\psi (u,v)\leq \frac{C}{1+|(u,v)|^{S+2d+1}} 
$. From \eqref{OLD_1A}, \eqref{OLD_1B} and the above computations we get
\begin{eqnarray}
|H_h(h,(x,y),(x',y'))|&\le& Ch^{-1/2} h^{-2d}\frac{1+|u|^2+|v|^2}{(1+| (u,v)| ^{S+2d+1}) }\nonumber\\
&\le& C\rho^{-1}\zeta_\rho(x'-x,y'-y-\frac{x+x'}2\rho^2).
\label{CTR_ONE_STEP}
\end{eqnarray}

For $1<j\le N$, we proceed like in the proof of Lemma 3.9
in \cite{kona:mamm:00}. We get that
\begin{equation*}
H_{h}(t_j,(x,y),(x',y'))=(\widehat H_{h}-\widetilde H_h)(t_j,(x,y),(x',y'))
\end{equation*}
where
\begin{eqnarray}
&&\widehat H_{h}(t_j,(x,y),(x',y'))=h^{-1}\int q_n\left( \theta_{1},\theta
_{2}\right)\times\nonumber\\
  &&\hspace*{30pt} \left[ \lambda (x+\widehat\gamma ^{1}(\theta_{1}),y+x h+ \widehat\gamma
^{2}\left( \theta_2\right) )-\lambda (x,y+x h)\right]
d\theta_{1}d\theta_{2},  \label{anc_14}
\end{eqnarray}
\begin{eqnarray}
&&\widetilde H_{h}(t_j,(x,y),(x^{\prime },y^{\prime
}))=h^{-1}\int q_n\left( \theta
_{1},\theta_{2}\right)\times  \nonumber\\
&&\hspace*{30pt}\left[
\lambda (x+\widetilde{\gamma }^{1}(\theta_{1}),y+xh+\widetilde{\gamma }%
^{2}\left( \theta_{2}\right) )-\lambda (x,y+xh)\right]
d\theta_{1}d\theta_{2},  \label{anc_15}
\end{eqnarray}
with $\lambda(u,v)=\widetilde p_h(t_{j}-h ,(u,v),(x',y')) $,
\begin{eqnarray*}
\widehat \gamma^1(\theta_1)
=
hb(x,y)+\sqrt h\sigma(x,y)\theta_1,\ 
 \widehat \gamma^2(\theta_2)= \left\{ \frac{b(x,y)\gamma_nh}{2}+\sqrt h\sigma(x,y)\theta_2 \right\}h,
\end{eqnarray*}
and
\begin{eqnarray*}
\widetilde \gamma^1(\theta_1)&=&hb(x',y')+\sqrt h\sigma(x',y'-x'\rho^2)\theta_1,\\
\widetilde \gamma^2(\theta_2)&=& \left\{
\frac{b(x',y')\gamma_nh}{2}+\sqrt h\sigma(x',y'-x'\rho^2)\theta_2 \right\}h.
\end{eqnarray*}

Using a Taylor expansion at order three for $\lambda $ in
\eqref{anc_14} and \eqref{anc_15} we obtain
\begin{eqnarray*}
&&H_{h}(t_j,(x,y),(x^{\prime },y^{\prime }))=\biggl\{\frac{1}{2}%
\Tr(a(x,y)-a(x',y'-x'\rho^2 ))D^2_x\lambda (x,y+xh)\\
&&\hspace*{-20pt}+\langle b(x,y)-b(x',y'),\nabla_x\lambda (x,y+xh)+\frac{\gamma_n h}{2}\nabla_y\lambda (x,y+xh)\rangle\biggr\}\\
&&+ 
\biggl\{ \frac{h}{2}\left( \langle
D^2_{x}\lambda(x,y+xh)b(x,y),b(x,y)\rangle -\langle
D^2_{x}\lambda(x,y+xh)b(x',y'),b(x',y')
\rangle\right)\\
&&+\frac{h^2}{2}\left(  \tr(  D^2_{y}\lambda(x,y+xh) (a(x,y)-a(x',y'-\rho^2x')) ) \right)\frac{1}{3}\gamma _{n}\left( 1+\frac{1}{2n}\right)\\
&&+\frac{\gamma_n^2 h^3}{8}\left(\langle D^2_y\lambda(x,y+xh)b(x,y),b(x,y) \rangle -\right.  
\left. \langle D^2_y\lambda(x,y+xh)b(x',y'),b(x',y')\rangle \right)
\biggr\}+\\
&   ²&\{h^{-1}\bint{}^{}d\theta_1d\theta_2d q_n(\theta_1,\theta_2)\\
&&\left(\langle D^2_{y,x}\lambda(x,y+xh)\widehat \gamma^1(\theta_1),\widehat \gamma^2(\theta_2) \rangle
-\langle D^2_{y,x}\lambda(x,y+xh)\widetilde
\gamma^1(\theta_1),\widetilde \gamma^2(\theta_2) \rangle
\right)\}
\end{eqnarray*}
\begin{eqnarray*}
+3 h^{-1}\sum_{| \nu | =3}\int d\theta_1d\theta_2  \int_{0}^{1}d\delta(1-\delta
)^{2}q_{n}( \theta_{1},\theta_{2})
\frac{(\widehat\gamma ^{1}(\theta_{1}),\widehat \gamma ^{2}(\theta_2)^{\nu
}}{\nu!}\\
D^{\nu }\lambda (x+\delta \widehat \gamma
^{1}(\theta_{1}),y+xh+\delta \widehat\gamma
^{2}(\theta_2))
\end{eqnarray*}
\begin{eqnarray}
-3h^{-1}\sum_{| \nu | =3}\int d\theta_1d\theta_2 \int_{0}^{1}d\delta(1-\delta
)^{2}q_n( \theta_{1},\theta_{2})
\frac{(\widetilde\gamma ^{1}(\theta_{1}),\widetilde \gamma ^{2}(\theta_2))^{\nu
}}{\nu!}\nonumber\\
\times D^{\nu }\lambda (x+\delta \widetilde \gamma
^{1}(\theta_{1}),y+xh+\delta \widetilde\gamma
^{2}(\theta_2))
\nonumber\\
:=I+II+III+IV-V,  \label{anc_16}
\end{eqnarray}
where we denote $D^2_x\lambda(x,y+xh) $ (resp. $D^2_y\lambda(x,y+xh),\
D^2_{y,x}\lambda(x,y+xh) $) the $\R^d\otimes\R^d $  matrices
$(\partial_{x_i,x_j}\lambda(x,y+xh))_{(i,j)\in\leftB 1,d \rightB^2}
$ (resp. $(\partial_{y_i,y_j}\lambda(x,y+xh))_{(i,j)\in\leftB 1,d
\rightB^2}$, $(\partial_{y_i,x_j}\lambda(x,y+xh))_{(i,j)\in\leftB
1,d\rightB^2 } $). 

In the sequel, a useful result is the following.
There exists $C>0$ s.t. for multi-indices $\alpha,\beta,\ |\alpha|\le 3,|\beta|\le 3 $,
\begin{eqnarray}
| \partial_x^\alpha \partial_{y}^\beta\lambda (x,y+xh) |&\leq&
C\rho ^{-(|\alpha|+3|\beta|)}\zeta _{\rho }\left( x^{\prime
}-x,y^{\prime }-y-xh\right.\nonumber\\
&&\left.-\frac{x+x'}{2}(\rho^2-h)\right) \nonumber \\
&\leq& C\rho ^{-(|\alpha|+3|\beta|)}\zeta _{\rho }\left( x^{\prime
}-x,y^{\prime }-y-\frac{x+x'}{2}\rho^2\right) \label{ancien_17}.
\end{eqnarray}
This assertion can be proved similarly to Lemma 3.7 in \cite{kona:mamm:00}.

Note now that
\begin{eqnarray*}
I&=&(L-\widetilde{L})\widetilde
p_h(t_j,(x,y+xh),(x',y'))\\
&&+\biggl\{ \frac{h\gamma_n}{2}\langle b(x,y)-b(x',y')  
,\nabla_{y}\lambda(x,y+xh)\rangle
\\
&&+
(L-\widetilde{L} )\bigl(\lambda (x,y+xh)
-\widetilde p_h(t_j,(x,y+xh),(x',y')) \bigr)\biggr\}\\
&:=&(K_{h}
+M_h^1)(t_j,(x,y),(x^{\prime
},y^{\prime })).
\end{eqnarray*}
From the above equation and \eqref{ancien_17} we get
\begin{eqnarray}
\label{SECOND_BOUND}
| M^1_h(t_j,(x,y),(x',y'))|&\le& C
\rho^{-1}
\zeta_\rho(x'-x,y'-y-\frac{x+x'}{2}\rho^2).\nonumber\\
\end{eqnarray}

Using similarly \eqref{ancien_17} and tedious but elementary calculations, one can split in $II,III$ the terms that give the expected order, i.e. bounded by $C\sqrt h \rho^{-1}\zeta_\rho(x'-x,y'-y-\frac{x+x'}{2}\rho^2) $ and denoted below by $R_h^2(t_j,(x,y),(x',y')) $, and those that give an integrable singularity in time, i.e. bounded by $C\rho^{-1}\zeta_\rho(x'-x,y'-y-\frac{x+x'}{2}\rho^2) $ and denoted below by $M_h^2(t_j,(x,y),(x',y'))$.

It remains to estimate $IV-V$ in \eqref{anc_16}.  To this end
write,
\begin{eqnarray*}
 IV-V=3h^{-1}\sum_{| \nu  | =3}%
\frac{1}{\nu !}\int d\theta_1d\theta_2\int_{0}^{1}d\delta (1-\delta )^{2} q_n(\theta_1,\theta_2)\biggl\{\\
((\widetilde \gamma^1(\theta_1),\widetilde \gamma^2(\theta_2))^\nu-(\widehat \gamma^1(\theta_1),\widehat \gamma^2(\theta_2))^\nu)D^\nu\lambda(x+\delta \widetilde \gamma^1(\theta_1),y+xh+\delta \widetilde \gamma^2(\theta_2))\\
+(\widehat \gamma^1(\theta_1),\widehat
\gamma^2(\theta_2,))^\nu\bsum{|\mu|=1}^{}\bint{0}^{1}d\alpha
D^{\nu,\mu}\lambda(x+\delta \widehat \gamma^1(\theta_1)+\alpha\delta
(\widetilde \gamma^1-\widehat \gamma^1)(\theta_1)
,\\
y+xh+\delta \widehat \gamma^2(\theta_2)+\alpha\delta (\widetilde \gamma^2(\theta_2)-\widehat \gamma^2(\theta_2)) )\\
\left(\delta (\widetilde \gamma^1-\widehat \gamma^1)(\theta_1) ,\delta (\widetilde \gamma^2(\theta_2)-\widehat \gamma^2(\theta_2))\right)^\mu
\biggl\}:=M_h^3(t_{j},(x,y),(x',y')).
\end{eqnarray*}
Computations involving \eqref{ancien_17} yield
\begin{align}
|M_h^3(t_j,(x,y),(x',y'))|
 \leq C\rho ^{-1}\zeta _{\rho }(x^{\prime }-x,y^{\prime }-y-\frac{x+x'}{2}\rho^2).
\label{DEF_M_4_H}
\end{align}
We refer to the proof of (3.80) p. 584 in \cite{kona:mamm:00} and Appendix \ref{APP_TRICK} for additional details.
This completes the proof.\finpreuve

The proof of Lemma \ref{OLD_LEMME2} then follows from the previous
proof, \eqref{ancien_17}, \eqref{SECOND_BOUND}, \eqref{DEF_M_4_H}
and \eqref{anc_16} for $j\in (1,N\rightB $ and \eqref{CTR_ONE_STEP} for
$j=1 $.

\section{Control of the $(S_i)_{i\in\leftB 1,4 \rightB}$}
\label{CTR_SI}

\subsection{Control of $S_1$} Set $$ p_d(T,(x,y),(x',y'))=\sum_{r=0}^{\infty} \tilde p\otimes_h H^{(r)}(T,(x,y),(x',y')) .$$ From Proposition \ref{PARAM_ANNEXE} one has
\begin{eqnarray*}
( p-p_d)(T,(x,y),(x',y'))&=&( p\otimes H- p\otimes_h H)(T,(x,y),(x',y'))\\
&&+( p- p_d)\otimes_h H(T,(x,y),(x',y')).
\end{eqnarray*}
Iterating the previous identity we get
\begin{eqnarray}
( p- p_d)(T,(x,y),(x',y'))=( p\otimes H- p\otimes_h H ) \otimes_h \varphi(T,(x,y),(x',y')),\nonumber\\
\label{CTR_DIFF_INF}
\end{eqnarray}
where $\forall j\in\leftB 0,N-1\rightB,\ \forall (u,v)\in\R^{2d}$,
\begin{eqnarray*}
\label{DEF_VARPHI} 
\varphi(T-t_j,(u,v),(x',y'))=\bsum{r=0}^{\infty} H_h^{(r)}(T-t_j,(u,v),(x',y')).
\end{eqnarray*}
Let us first give a bound for $P_j(u,v):=( p\otimes H- p\otimes_h H
)(t_j,(x,y),(u,v)),\ j\in\leftB 0,N\rightB, \
(u,v)\in\R^{2d}$. First, from the previous definitions of the
continuous and discrete convolution operators, $P_0(u,v)=0 $, in the
sense of generalized functions. For $j\ge 1 $ write
\begin{eqnarray*}
P_j(u,v)&=&\bsum{i=0}^{j-1}\bint{t_i}^{t_{i+1}}dt\bint{\R^{2d} }^{}dw dz \lambda_{(u,v)}(t,(w,z))-\lambda_{(u,v)}(t_i,(w,z)),\\
\lambda_{(u,v)}(t,(w,z))&:=&  p(t,(x,y),(w,z))H(t_j-t,(w,z),(u,v)). 
\end{eqnarray*}
A first order Taylor expansion and Fubini's theorem give
\begin{eqnarray}
P_j(u,v)&=&\bsum{i=1}^{j-1}\bint{t_i}^{t_{i+1}}dt (t-t_i)\bint{0}^{1}d\delta Q_i^\delta(u,v,s)
        +T_j^0,\nonumber\\
Q_i^\delta(u,v,s)&:=&\bint{\R^{2d}}^{}dw dz  \partial_s\lambda_{(u,v)}(s,(w,z))_{s=t_i+\delta(t-t_i)},\ i\in\leftB 1,j-1\rightB.\nonumber\\
T_j^0&:=& \bint{0}^h dt \bint{\R^{2d} }^{}dwdz   p(t,(x,y),(w,z))\nonumber\\
&&\times(H(t_j-t,(w,z),(u,v))-H(t_j,(x,y),(u,v))).
\label{EXP_PUV}
\end{eqnarray}
From Lemma \ref{CTR_DENS_L}, Theorem \ref{THEO_DIFF}, we obtain
\begin{eqnarray*}
T_j^0&\le & C\sqrt h \widehat p_c( t_j,(x,y),(u,v)),\\
|\partial_s\lambda_{(u,v)}(s,w,z)  |&\le&
C(s^{-1}(t_j-s)^{-1/2}+(t_j-s)^{-3/2})\\
&&\times \widehat p_c(s,(x,y),(w,z))\widehat p_c( t_j-s,(w,z),(u,v)).
\end{eqnarray*}
Plug now the above control in \eqref{EXP_PUV}, we get
\begin{eqnarray*}
P_j(u,v)&\le& C   
\widehat p_c(t_j,(x,y),(u,v))(h^{1/2}+h^{2}\left(t_j^{-1/2}\bsum{i=1}^{\lfloor (j-1)/2\rfloor}t_i^{-1}\right.\\
&& \left. +t_j^{-1}\bsum{i=\lfloor (j-1)/2\rfloor+1}^{j-2}(t_j-t_{i+1})^{-1/2}+\bsum{i=1}^{j-2}t_i^{-3/2}\right) )\\
&\le & Ch^{1/2}  
\widehat p_c(t_j,(x,y),(u,v)).
\end{eqnarray*}
Hence, from \eqref{CTR_DIFF_INF} and a suitable version of \eqref{THE_EQ_POUR_LES_QUEUES} for the discrete convolution operator we derive
\begin{eqnarray*}
|( p- p_d)(T,(x,y),(x',y'))|&\le & Ch^{1/2}  
\widehat p_c(T,(x,y),(x',y')).
\end{eqnarray*}
The bound for $S_1$ can be derived using once again \eqref{THE_EQ_POUR_LES_QUEUES} for both the continuous and discrete convolution operators and the asymptotics of the Gamma function.\\

\subsection{Control of $S_2$}
Define for $r\in\leftB 0,N\rightB,\ T_r:=(\widetilde p \otimes H^{(r)}-\widetilde p \otimes \widetilde H^{(r)} )(T,(x,y),(x',y'))$. For $r=1$, with the notations of Lemma \ref{CTR_DENS_L} one gets
\begin{eqnarray*}
|T_1|&\le & Ch^2\bsum{j=1}^{N-1} \bint{0}^{1}d\delta\bint{}^{}t_j^{-3/2}(T-t_j)^{-1/2}|u|\widehat p_c(t_j,(x,y),(u,v))\\&& \times \widehat p_c(T-t_j,(u,v+\delta u h),(x',y'))dudv.
\end{eqnarray*} 
Also, for a different constant $c'$ than the one appearing in $\widehat p_c $, $\forall j\in\leftB 1,N-1\rightB $,
\begin{eqnarray}
&&\widehat p_c(T-t_j,(u,v+\delta u h),(x',y')) \le C(T-t_j)^{-(3k+d)/2}\nonumber \\
&&\times \exp\left(-c'\left\{\frac{|x'-u|^2}{4(T-t_j)} +3\frac{|y'-v-\frac{u+x'}{2}(T+\delta h-t_j)|^2}{(T-t_j)^3}\right\} \right)\nonumber \\
&\le& C\widehat p_{c'}( T+\delta h-t_j,(u,v),(x',y')). \label{TRICK_TIME}
\end{eqnarray}
Hence, up to another suitable modification of the constant in order to have the semigroup property
\begin{eqnarray*}
|T_1|&\le & Ch^2\bsum{j=1}^{N-1} t_j^{-3/2}(1+(T-t_j)^{-1/2}|x'|)\bint{0}^{1}d\delta \widehat p_c(T+\delta h,(x,y),(x',y'))\\
&\le & Ch^{1/2}(1+T^{-1/2}|x'|)\bint{0}^{1}d\delta \widehat p_c(T+\delta h,(x,y),(x',y')).
\end{eqnarray*}
Write now, for all $r\ge 2 $,
\begin{eqnarray*}
T_r&=&\widetilde p\otimes_h H^{(r-1)} \otimes_h (H-\widetilde H)(T,(x,y),(x',y'))\\
&& +(\widetilde p\otimes_h H^{(r-1)}-\widetilde p\otimes_h \widetilde H^{(r-1)})\otimes_h \widetilde H(T,(x,y),(x',y')):=T_{r1}+T_{r2}.
\end{eqnarray*}
The term $T_{r1} $ can be handled as $T_1$ exploiting the control 
\begin{eqnarray*}
|\partial_z \widetilde p\otimes_h H^{(r)}(t,(x,y),(w,z))|&\le& C^{r+1}t^{(r-3)/2}\widehat p_c(t,(x,y),(w,z))\\
&&\times \prod_{i=0}^{r}B((1+i)/2,1/2).
\end{eqnarray*}
For $T_{r2}$ one uses the control of step $(r-1)$. Completing the induction one derives
$$|S_2|\le Ch^{1/2}(1+|x'|)\sup_{\delta\in[0,1]}\widehat p_c(T(1+\delta),(x,y),(x',y')) .$$
Note that this term is the only one for which we have a linear contribution of the terminal variable. This is, because of the shift, in some sense unavoidable. Also, the previous trick
in \eqref{TRICK_TIME} adds the constraint to take a supremum w.r.t. to a twice larger time interval
as the initial one.

\subsection{Control of $S_3$}
\label{APP_TRICK}
For $r=1$ we have to control
\begin{eqnarray*}
\widetilde p\otimes_h M_h(T,(x,y),(x',y'))=\\
\bsum{i=1}^{3}h\bsum{j=0}^{N-1}\bint{}^{}dudv\widetilde p(t_j,(x,y),(u,v))M_h^i( T-t_j,(u,v),(x',y'))\\
:=h\bsum{i=1}^{3}\bsum{j=0}^{N-2}I_{i,j}+hI_{N-1}.
\end{eqnarray*}

 The term $hI_{N-1}$ needs to be handled by a different technique than the other ones. Write
\begin{eqnarray*}
hI_{N-1}&=&\bint{}^{}du dv \widetilde p(T-h,(x,y),(u,v))\\
&&\times\left(p_h
-\widetilde p_h\right)(h,(u,v),(x',y')).
\end{eqnarray*}
Set $V=(\frac{u-x}{(T-h)^{1/2}},\frac{v-y-\frac{x+u}{2}(T-h)}{(T-h)^{3/2}})$. Write now $|u-x|=|x'-x+u-x'| $,
\begin{eqnarray*}
\left|v-y-\frac{x+u}{2}(T-h)\right|=\left|y'-y-\frac{x+x'}2 T+ \frac{x-x'}2h +v-y'+uh\right.\\
\left.+\frac{x'-u}2(T+h)\right|
\end{eqnarray*}
Set  $U=\left(\frac{x'-x}{(T-h)^{1/2}},\frac{y'-y-\frac{x+x'}{2}T+\frac{x-x'}{2}h}{(T-h)^{3/2}}\right), V:=U+R $.
Recall also from Lemma \ref{CTR_DENS_L} that $\widetilde p\le C \widehat p_c$. Hence, for all $Z\in \N^*,\ \exists C:=C(Z),\ \widetilde p(T-h,(x,y),(u,v))\le (T-h)^{-2d}\frac{C}{1+|V|^Z}$.
From the basic identity $\frac{1}{1+|U+R|^Z}\le \frac{\max (2^Z,1+(2|R|)^Z)}{1+|U|^Z} $ and the definitions of the models \eqref{CDM} and \eqref{CDM_FRO}, using the same techniques as in the proof of Lemma \ref{CTR_KERNL_MARK} for the case $j=1$ one gets:
\begin{eqnarray*}
|hI_{N-1}|&\le& \frac{Ch^{1/2}}{1+|U|^Z}\bint{}^{}du'dv'(1+ |(u',v')|^{Z})\psi(-u',-v')\\
&\le& Ch^{1/2}\zeta_{\sqrt
{T}}(x'-x,y'-y-\frac{x+x'}{2}T),\\
\end{eqnarray*}
taking $Z=S-4 $ for the last inequality. 

Also, from the definitions of the $(M_h^i)_{i\in\leftB 1,3\rightB} $
in the previous section and using freely its notations, we derive
for all $j\in\leftB 0,N-2\rightB $:
\begin{eqnarray*}
|M_h^1( T-t_j,(u,v),(x',y'))|&\le&
h(T-t_j)^{-3/2}\zeta_\rho(x'-u,y'-v-\frac{(u+x')}{2}(T-t_j)),
\end{eqnarray*}
from which one gets $h\bsum{j=0}^{N-2}|I_{1,j}|\le Ch^{1/2}\zeta_{\sqrt T}(x'-x,y'-y-\frac{x+x'}{2}T) $. The terms in $M_h^2$ coming from $II$ in \eqref{anc_16} can be handled as $M_h^1$. For those coming from $III$, i.e. crossed derivatives, the contribution associated to $j=0$ is easily analyzed and for $j>1$ an integration by part w.r.t. $u$ leads to the same control. 
The trickiest term to analyze is $M_h^3$. Exploiting thoroughly
\eqref{ancien_17} and Lemma \ref{CTR_DENS_L}, the proof is similar
to the one in \cite{kona:mamm:00}, see p.578 control of (3.45), that
relies on suitable integration by parts. We omit the details here.
Actually, for $r\ge 1$ it can be shown by induction that
\begin{eqnarray*}
\left|\left(\widetilde p\otimes_h \widetilde H^{(r)}-\widetilde p\otimes_h(\widetilde H+M)^{(r)}\right)(T,(x,y),(x',y')) \right|\\
\le  \frac{h^{1/2}C^{r+1}}{\Gamma([r+2]/2)}\chi_{\sqrt
T}(x'-x,y'-y-\frac{x+x'}{2}T),
\end{eqnarray*}
which gives the control.

\subsection{Control of $S_4$} One can show that Lemma \ref{THE_LEMME_DENS_GELEE} is still valid for the derivatives of the frozen densities. Using this result and Lemma \ref{OLD_LEMME2}, the proof is then similar to the one of \cite{kona:mamm:00}.

\bibliographystyle{alpha}
\bibliography{bibli}

\end{document}